\documentclass[journal, 12pt, onecolumn]{IEEEtran}

\usepackage{url}
\usepackage{blindtext}
\usepackage{subcaption}
\usepackage{multirow}
\usepackage{graphicx}
\usepackage{cite}
\usepackage{amsfonts}

\usepackage{graphicx}
\usepackage{amsmath,amssymb,amsfonts}
\usepackage{amsthm}
\usepackage{mathrsfs}
\usepackage{xcolor}
\usepackage{textcomp}
\usepackage{manyfoot}
\usepackage{listings}

\usepackage{subcaption}
\usepackage{amsmath}

\usepackage{multicol}
\usepackage{multirow}
\usepackage{algpseudocode}
\usepackage{algorithm}

\usepackage{caption}
\usepackage{subcaption}
\usepackage{amsmath}
\usepackage{multicol}
\usepackage{xcolor,soul,framed} 
\usepackage{amssymb}
\usepackage{algorithm}
\usepackage{algpseudocode}
\usepackage{graphicx}
\usepackage{multirow}
\usepackage{tabularx} 
\usepackage{lscape}   
\usepackage{xcolor}
\definecolor{red}{RGB}{255,0,0}

\begin{document}

\title{Accelerating genetic optimization of nonlinear model predictive control by learning optimal search space size}

\author{Eslam~Mostafa,
        Hussein~A.~Aly,
        Ahmed~Elliethy
\thanks{Authors are with the Department of Electrical and Computer Engineering, Military~Technical~College, Cairo 14627, Egypt (e-mails: eslammostafafawzy@mtc.eg.edu, haly@ieee.org, a.s.elliethy@mtc.eg.edu).}
\thanks{This    paper    has    supplementary    downloadable    material    available  provided  by  the  authors.  The  material  consists  of a supplementary document.}
\thanks{Color  versions  of  one  or  more  of  the  figures  in  this  paper  are  available in the electronic version of this manuscript.}}


\maketitle

\begin{abstract}

Genetic algorithm (GA) is typically used to solve nonlinear model predictive control's optimization problem. However, the size of the search space in which the GA searches for the optimal control inputs is crucial for its applicability to fast-response systems. This paper proposes accelerating the genetic optimization of NMPC by learning optimal search space size. The approach trains a multivariate regression model to adaptively predict the best smallest size of the search space in every control cycle. The proposed approach reduces the GA's computational time, improves the chance of convergence to better control inputs, and provides a stable and feasible solution. The proposed approach was evaluated on three nonlinear systems and compared to four other evolutionary algorithms implemented in a processor-in-the-loop fashion. The results show that the proposed approach provides a 17-45\% reduction in computational time and increases the convergence rate by 35-47\%. The source code is available on GitHub.

\end{abstract}


\newcommand{\grx}{\chi}
\newcommand{\gry}{\zeta}
\newcommand{\rgrx}{^r\grx}
\newcommand{\rgry}{^r\gry}
\newcommand{\betaB}{{\boldsymbol\beta}}
\newcommand{\omegaB}{{\boldsymbol\omega}}
\newcommand{\gammaB}{{\boldsymbol\gamma}}
\newcommand{\alphaB}{{\boldsymbol\alpha}}
\newcommand{\states}{\mathbf{x}}
\newcommand{\refx}{\mathbf{r}}
\newcommand{\refxi}{r}
\newcommand{\refu}{\mathbf{v}}
\newcommand{\refui}{v}
\newcommand{\optIns}{\mathbf{z}_c^*}
\newcommand{\ins}{\mathbf{z}_c}
\newcommand{\hzn}{h}
\newcommand{\bsm}{\delta_c}
\newcommand{\bsmm}{\delta}
\newcommand{\BSM}{\boldsymbol\Delta_{c}}
\newcommand{\BSMmax}{\Delta_{c}^{\max}}
\newcommand{\Cmax}{\mathcal{C}_{c}^{\max}}
\newcommand{\Emax}{\mathcal{E}_{c}^{\max}}
\newcommand{\Cc}{\mathbf{\mathcal{C}}_{c}}
\newcommand{\Ec}{\mathbf{\mathcal{E}}_{c}}
\newcommand{\E}{\mathbf{\mathcal{E}}}
\newcommand{\Wc}{\mathbf{\mathcal{W}}_{c}}
\newcommand{\Rc}{\mathbf{\mathcal{R}}_{c}}
\newcommand{\SampDens}{\nu}
\newcommand{\kernelScale}{\gamma}

\newcommand{\finalBSM}{\boldsymbol\Psi{c}}

\newcommand{\ArxivFig}{./ArxivFig}
%
\IEEEpeerreviewmaketitle
\section{Introduction}\label{sec1}

Model predictive control (MPC) is a robust method for controlling multivariate systems while satisfying constraints~\cite{garcia1989MPCbook}. It generates optimal control inputs by solving a multivariate optimization problem at each control cycle, considering both current and future states. The key strength of MPC is its ability to optimize control actions over a future time horizon. However, the computational complexity of MPC increases significantly for complex or nonlinear systems due to nonlinearities, coupled inputs, and safety constraints, making the optimization non-convex. This challenge is especially evident when MPC runs on resource-constrained embedded hardware or is required to control systems requiring fast response.

Various classical optimization solvers, both linear and nonlinear, have been developed to address these challenges~\cite{gelleschus2018solverscomparison}. While linear solvers assume system linearity, many systems require nonlinear solvers to handle their nonlinear dynamics and constraints. Unfortunately, this leads to higher computational times, which can be problematic for systems with short control cycles. Other than these classical solvers, several evolutionary algorithms have been used to address optimization challenges in Nonlinear Model Predictive Control (NMPC) across various applications such as Particle Swarm Optimization (PSO)~\cite{jain2022PSOoverview}, Differential Evolution (DE)~\cite{deng2021DEoverview}, and the Genetic Algorithm (GA)~\cite{kramer2017geneticBook}. These evolutionary algorithms offer two key advantages compared with the classical solvers. First, they can handle non-differentiable, discontinuous, or nonlinear components without relying on gradient evaluation, making them versatile for various problems. Additionally, their computations can be accelerated through parallelization on GPUs. These advantages make these algorithms represent powerful tools for solving complex optimization problems.

The GA is a widespread evolutionary algorithm that is extensively used in NMPC~\cite{popov2005GAopt1,sivanandam2008GAopt2,karakativc2021GAopt3,deb2001GAopt4}. It solves optimization problems by evaluating a population of candidate solutions from a defined search space. Solutions evolve until the control cycle time expires or a termination condition is met. However, the search space size is crucial for systems requiring fast response. A small search space reduces the likelihood of finding optimal\footnote{Please note that global optimality is not guaranteed by the genetic algorithm when used to solve non-convex optimization problems in general. So, we use the term ``best'' instead of ``optimal'' control inputs. What we meant by ``best control inputs'' here are the control inputs that minimize the optimization cost function below a certain threshold.} control inputs, while a large space increases this probability but requires more computations~\cite{sarker2003popVSsearchSpace2,hassanat2019popVSsearchSpace2,bayas2021popVSsearchSpace1}. Carefully limiting the search space can help balance computational time and solution quality.

This paper presents an approach to accelerate the genetic optimization of nonlinear model predictive control by learning the optimal search space size. Our method trains a regression model to adaptively predict the optimal minimal search space size for each control cycle, improving the likelihood of finding the best control inputs within minimal computational time. To achieve this, we first generate a synthetic dataset using the system's model, considering the inaccuracies in the
actions applied to the system, model uncertainties, external disturbances, and changes in system dynamics. A regression model is then trained on this dataset to estimate the optimal smallest search space size. This estimated size is fed into the GA, guiding the search for control inputs within the defined space.

To demonstrate the practical effectiveness of our approach, we implemented it on an Nvidia\textsuperscript{\texttrademark} Jetson TX2~\cite{NvidiaSpec} embedded platform using a processor-in-the-loop (PIL) setup. The approach was evaluated on three nonlinear systems and compared to four state-of-the-art evolutionary algorithms for solving NMPC optimization, all implemented on the same platform. Our method achieved a 17-45\% reduction in computational time compared with the next best-compared algorithm while significantly increasing the likelihood of obtaining the best control inputs. Specifically, the convergence rate to the best control inputs before the termination of the cycle improved by 35-47\% compared to the next best algorithm. The source code for our approach is available online.

The remainder of this paper is organized as follows. Section~\ref{sec:survey} presents a short survey of using the GA in addition to several other evolutionary algorithms used to tackle the NMPC optimization in various application domains. Section~\ref{sec:nmpc_formulation_ga} presents the mathematical formulation of the NMPC and the GA algorithm in detail. In Section~\ref{sec:proposed}, we motivate and present the proposed approach for learning the optimal search space size for genetic optimization. Section~\ref{sec:results} presents the experimental setup and discusses the experimental results that evaluate our proposed approach. Section~\ref{sec:conclusion} summarizes our conclusion and presents our future works.

\section{Literature survey}
\label{sec:survey}

Numerous studies have explored the application of evolutionary algorithms in optimizing Nonlinear Model Predictive Control (NMPC), particularly for handling nonlinear and complex systems. This section provides an overview of these studies and highlights the integration of various evolutionary algorithms with NMPC.

Particle Swarm Optimization (PSO)~\cite{jain2022PSOoverview,chai2018PSO1}, inspired by the social behavior of animals, is commonly used to address the challenges of the optimization problem in NMPC. PSO enables efficient exploration of large search spaces, adjusting each particle's position based on individual and neighbor experiences. For example, in~\cite{chai2020PSO1Plann}, the PSO is integrated into NMPC to optimize control inputs for nonlinear systems in real-time, significantly reducing computational demands and enhancing stability and convergence. Similarly, the method in~\cite{chai2020Pso2Plann} applied a hybrid PSO for autonomous vehicle parking, balancing exploration and exploitation for more accurate optimization. Additionally, in~\cite{chai2019PSO2}, PSO's effectiveness is demonstrated in NMPC for vehicle trajectory planning, showcasing its ability to solve time-sensitive, multi-objective problems with high accuracy and reduced computational complexity.

Another effective approach for NMPC optimization is the Differential Evolution (DE) method~\cite{deng2021DEoverview}, which searches for the global optimum in highly nonlinear systems. DE evolves potential solutions through iterative operations. In~\cite{zhang2021DEcompare}, DE was applied to optimize control inputs for a single-link flexible-joint (FJ) robot in NMPC. Additionally,~\cite{deng2021DE2} introduced an adaptive DE algorithm that adjusts its parameters dynamically, improving search efficiency for local optimization problems within each time horizon.

Genetic Algorithms (GA) have also been widely applied to NMPC optimization across various domains, with two primary research directions. The first involves using GA to directly solve NMPC optimization problems, either to find optimal or sub-optimal solutions. Optimal approaches aim to find the best solution within the control cycle, while sub-optimal approaches stop once a satisfactory solution is found. For example, the method in~\cite{arrigoni2022trjplanningGAMPC} used GA to solve NMPC for autonomous vehicle trajectory planning, while in~\cite{du2016develGAsaftyEqn}, a cost function is introduced for balancing safety and comfort. Other works like~\cite{samsam2023fuel1} applied GA for trajectory optimization in space missions, while in~\cite{hyatt2017real} proposed a parallelized GA for real-time NMPC using GPUs. In cases where a fast response is needed, sub-optimal solutions are used. For instance, in~\cite{chen2009suboptimal1}, a method is developed for a stirred tank reactor where optimization stops when the cost function improves, while in~\cite{sharma2013suboptimal2} a similar approach is applied for USV autopilot design, ensuring feasible control signals within each cycle.

The second direction uses GA to tune NMPC parameters like prediction horizons and weighting factors, improving performance without trial and error. For instance, in~\cite{picotti2022tune1}, an NMPC controller is tuned for a virtual motorcycle, achieving better tracking accuracy and control effort. Similarly, in~\cite{yu2022tune2}, weight parameters are optimized for ship trajectory tracking, and in~\cite{yasini2023tune3}, weighting parameters are optimized to improve satellite motion control.

Practically, the discussed approaches reduce GA's computational complexity but often at the cost of controller performance. This reduction is typically achieved by either using shorter horizons for NMPC optimization or limiting iterations for sub-optimal solutions when optimal convergence is impossible within the control cycle time. 

This paper proposes a novel method that reduces GA's computational complexity while preserving NMPC performance. Our approach achieves this by learning the optimal search space size by adaptively predicting the best smallest search space in each control cycle. This improves the likelihood of finding better control inputs within a shorter computation time without sacrificing performance. While it doesn't guarantee optimal inputs (similar to the original GA), it increases the chances of achieving better results compared to traditional GA-based NMPC optimization.

While this study focuses on the GA, the proposed approach could be extended to other evolutionary optimization algorithms. However, we leave this direction to future work.

\section{Non linear model predictive control formulation and genetic algorithm}
\label{sec:nmpc_formulation_ga}
This section introduces the NMPC and discusses its general mathematical formulation. Then, we discuss the design of the GA and how it can be used to solve complicated NMPC optimization problems. We gather all abbreviations in Table \ref{tab:abbreviations} for better readability of the paper.

\begin{table}[t!]
\centering
\caption{Table of Abbreviations}
\begin{tabular}{|c|l|}
\hline
\textbf{Abbreviation} & \textbf{Definition} \\ \hline
NMPC   & Nonlinear Model Predictive Control \\
GA     & Genetic Algorithm \\
OG     & Original Genetic Algorithm \\
MG     & Modified Genetic Algorithm \\
PSO    & Particle Swarm Optimization \\
DE     & Differential Evolution \\
MNSVR  & Multivariate Nonlinear Support Vector Regressor \\
BSM    & Best Smallest Margin \\
$\Delta_c$ & Best smallest margin at control cycle $c$ \\
$\alpha_c$ & Maximum of normalized margins \\
$J$    & Cost function \\
$\mathcal{L}$    & Loss function \\
$V$    & Terminal cost \\
$X_f$  & Terminal set of allowed state values \\
$\mathbf{u}_{k}$  & Vector of control inputs at time step $k$ \\
$\mathbf{x}_{k}$  & Vector of system states at time step $k$ \\
$\kappa_f(x)$ & Control law \\
$\rho$ & Factor to determine the strength of control input noise \\
$\theta$ & Factor to determine the strength of measurement noise \\
$\zeta$ & Slack variable for the support vector regressor \\
$\kernelScale$ & Scaling factor of the Gaussian kernel   \\
$\SampDens$ & Maximum population size \\
$\xi$ & Minimum population size \\
\hline
\end{tabular}
\label{tab:abbreviations}
\end{table}

\subsection{Non linear model predictive control}
\label{sec:NMPC_formulation}
We consider the class of discrete-time nonlinear systems with the following general formulation
\begin{equation}
\mathbf{x}_{k+1}= f(\mathbf{x}_{k}, \mathbf{u}_{k}),
\label{eqn:statespaceform}
\end{equation}
where $k \in \mathbb{Z}$ is a discrete-time instant, $\mathbf{x}_{k} = [ x^1_{k}, \dots, x^m_{k} ] ^T$ represents the vector of states at time instant $k$, $\mathbf{u}_{k} = [ u^1_{k}, \dots, u^n_{k} ] ^T$ represents the vector of system inputs at $k$.

\begin{figure}[t]
	\centerline{\includegraphics[width=.7\textwidth]{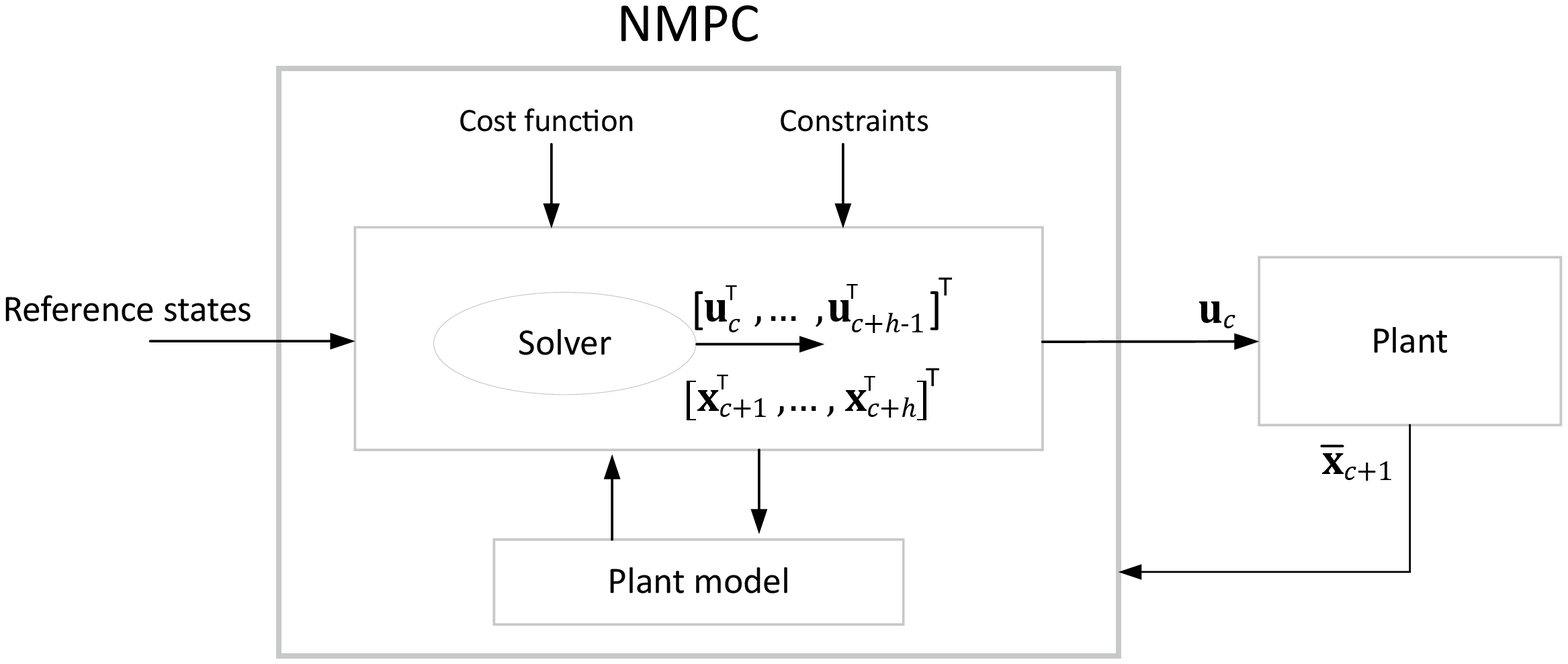}}
	\caption[Nonlinear model predictive control (NMPC) overall architecture.]{Nonlinear model predictive control (NMPC) overall architecture. At each control cycle, NMPC estimates the optimal vector of system inputs $ [\mathbf{u}^{T}_{c}, \dots, \mathbf{u}^{T}_{c + \hzn - 1}]^{T}$ that satisfies a given set of constraints for a fixed horizon length of $\hzn$ time steps. Only the first control input ($\mathbf{u}_{c}$) is applied to the system, and new system states ($\mathbf{\bar{x}}_{c+1}$) are produced, and this process is repeated in the next control cycle of NMPC.}
	\label{fig:NMPC}
\end{figure}

In each control cycle of the above system, the NMPC estimates the optimal vector of system inputs that minimizes a cost function $J$ and satisfies a given set of constraints over a fixed future horizon of length $\hzn$ time steps to reduce the deviation between the predicted states and the reference states $\refx_{k} = [\refxi^1_k, \dots, \refxi^{m}_k] ^{T}$ while considering the control efforts. As depicted in~Fig.~\ref{fig:NMPC}, in the $c^{\text{th}}$ control cycle, the NMPC estimates the system inputs $\ins = [\mathbf{u}^{T}_{c}, \dots, \mathbf{u}^{T}_{c + \hzn - 1}]^{T}$ over the future horizon. The NMPC estiamtes the inputs by employing the system model~\eqref{eqn:statespaceform} to generate the states $\mathbf{x}_{c+1}, \dots, \mathbf{x}_{c+h}$ over the future horizon from the current measured state $\mathbf{\bar{x}}_{c}$. Then, only the first control input $\mathbf{u}_{c}$ is applied to the actual plant system, and the new system states $\mathbf{\bar{x}}_{c+1}$ are measured. This process is repeated in all subsequent control cycles of the NMPC. 

The cost function $J$ comprises two terms: a loss function $\mathcal{L}$ and a terminal cost $V$. The function $\mathcal{L}$ penalizes the error between the current states with respect to their reference states and penalizes the control inputs within the optimization horizon. The terminal cost is chosen for stability insurance\cite{mayne2000stabilityGenal} and is not directly related to performance specification. Mathematically, the NMPC formulates and solves the following constrained optimization problem to estimate the optimal vector of system inputs $\optIns$. Specifically, 
\begin{equation}
\begin{aligned}
&\optIns = \arg \min_{\ins} J,\\
&\text { with } \\
&J := \sum_{k=c}^{c+\hzn-1} \mathcal{L}(\mathbf{x}_{k},\mathbf{u}_{k}) + V(\bold{x}_{c+h}),\\
&\text { s.t. } \\
&\mathbf{x}_{k+1}= f(\mathbf{x}_{k}, \mathbf{u}_{k}),\\
&\bold{x}_{c+h} \in X_f,\\
&\mathbf{x}_{k} \in \mathcal{X} \; \; \forall k,\\
&\mathbf{u}_{k} \in \mathcal{U} \; \; \forall k,\\
\end{aligned}
\label{eqn:MPCCost}
\end{equation}
where the constraints
\begin{equation}
\begin{aligned}
&\mathcal{X} = \{ \mathbf{b} :  \mathbf{x}_{\text{min}} \leq \mathbf{b} \leq \mathbf{x}_{\text{max}} \}, \\
&\mathcal{U} = \{ \mathbf{a} :  \mathbf{u}_{\text{min}} \leq \mathbf{a} \leq \mathbf{u}_{\text{max}} \}, \\
\end{aligned}
\end{equation}
are set to limit the states and the control inputs, respectively, to certain ranges, and $X_f$ represents the terminal set, which contains the allowed state values that the control aims to achieve at the end of the prediction horizon.

Note that the terminal cost $V$ and the terminal region $X_f$ are determined offline to enforce stability. Specifically, $V$ and $X_f$ are defined in such a way as to satisfy the following conditions to ensure stability.
\begin{itemize}
\item A1: State Constraint Satisfaction
\[
X_f \subseteq \mathcal{X}, \, X_f \text{ closed}, \, 0 \in X_f.
\]
This condition ensures that the terminal set \( X_f \) is a subset of the state constraint set \( \mathcal{X} \), is closed, and includes the origin. This guarantees that the state constraints are satisfied within the terminal set.

\item A2: Control Constraint Satisfaction
\[
\kappa_f(\mathbf{x}) \in \mathcal{U}, \, \forall \mathbf{x} \in X_f
\]
This condition ensures that the control action \( \kappa_f(\mathbf{x}) \) applied within the terminal set \( X_f \) respects the control constraints \( \mathcal{U} \).

\item A3: Positive Invariance
\[
f(x, \kappa_f(\mathbf{x})) \in X_f, \, \forall \mathbf{x} \in X_f
\]
This condition ensures that the terminal set \( X_f \) is positively invariant under the control law \( \kappa_f(\mathbf{x}) \). This means that if the system state starts within \( X_f \), it will remain within \( X_f \) under the control law \( \kappa_f(\mathbf{x}) \).

\item A4: Lyapunov Function Decrease
\[
V(f(\mathbf{x}, \kappa_f(\mathbf{x})))-V(\mathbf{x}) + \mathcal{L}(\mathbf{x}, \kappa_f(\mathbf{x})) \leq 0, \, \forall \mathbf{x} \in X_f
\]
This condition ensures that the Lyapunov function decreases over time within the terminal set \( X_f \), ensuring the stability of the closed-loop system.

\end{itemize}

The subsequent section presents the genetic algorithm, demonstrating the potential to solve complex optimization problems.

\subsection{Genetic algorithm as a NMPC solver}
\label{sec:ga-npmc-solver}

\begin{figure*}[t]
	\centerline{\includegraphics[width=.8\textwidth]{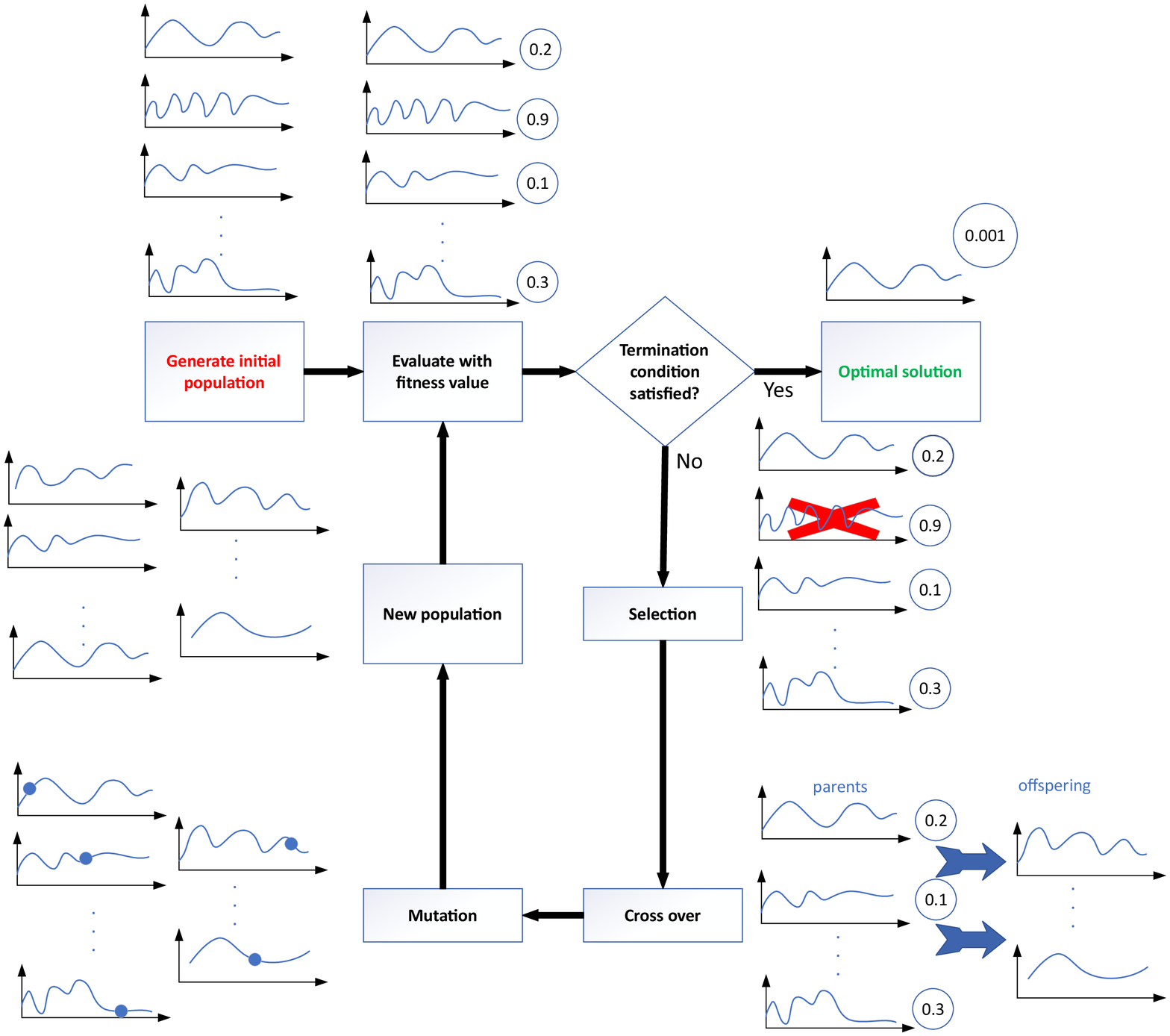}}
	\caption{Illustration of using the genetic algorithm to solve the NMPC optimization problem. GA begins by initializing a population of candidate solutions, which are evaluated, and the best are selected. These selected solutions then undergo crossover and mutation to produce a new generation. The process repeats until a termination condition or convergence is reached.}
	\label{fig:generalGenetic}
\end{figure*}

The GA is highly effective for solving complex NMPC optimization problems, particularly those involving nonlinearities or if/else conditions, where classical solvers struggle. A general GA architecture is illustrated in Fig.~\ref{fig:generalGenetic}. The algorithm begins with a population of candidate solutions (chromosomes), where each solution $\ins$ consists of $\hzn n$ values (genes). Solutions evolve through generations, with new offspring produced via crossover (combining solutions) and mutation (randomly altering genes). The crossover and mutation rates control how frequently these changes are made. After each generation, solutions are ranked by their fitness, determined by the cost function $J$ as:
\begin{equation}
    F = \frac{1}{(1+J)}.
    \label{eq:fitness}
\end{equation}
Selection methods, like a roulette wheel or tournament selection, favor higher fitness solutions for survival in future generations~\cite{jebari2013GAselection}. The process continues until a termination condition is met, either when the solution's cost falls below a threshold or when the allowed control cycle time runs out, returning the best solution found.

The GA's exploratory nature can be time-consuming, particularly when dealing with large solution spaces. This limits its use in applications needing fast responses. In the next section, we introduce our proposed approach to accelerate GA by adaptively reducing the search space size, enabling it to find better solutions faster.


\section{Proposed learning of optimal search space size for genetic optimization}
\label{sec:proposed}

\begin{figure*}[t]
\centering
\begin{subfigure}[]{0.48\textwidth}
   \centerline{\includegraphics[width=.7\textwidth]{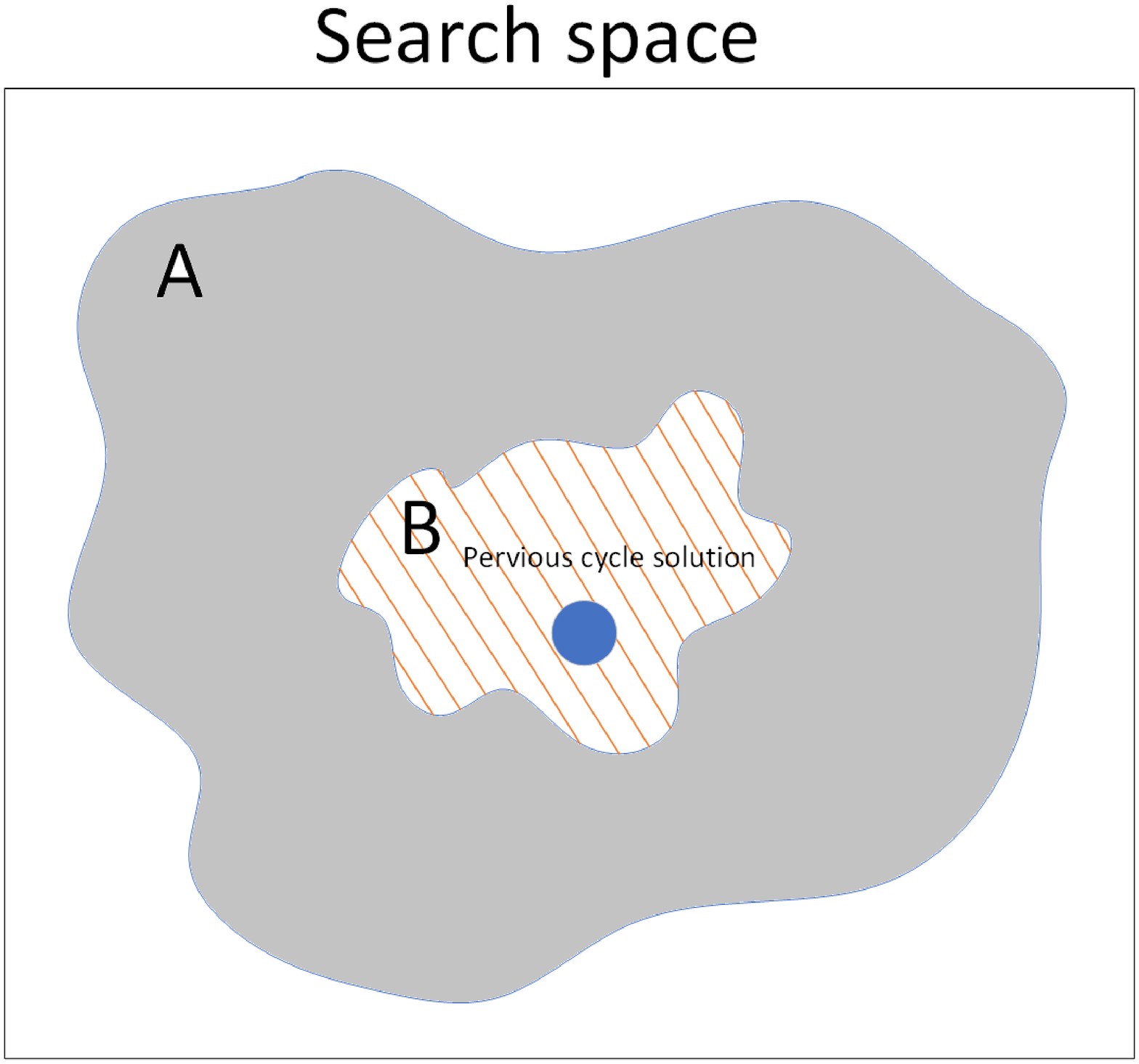}}
   \caption{}
\end{subfigure}
\begin{subfigure}[]{0.48\textwidth}
   \centerline{\includegraphics[width=\textwidth]{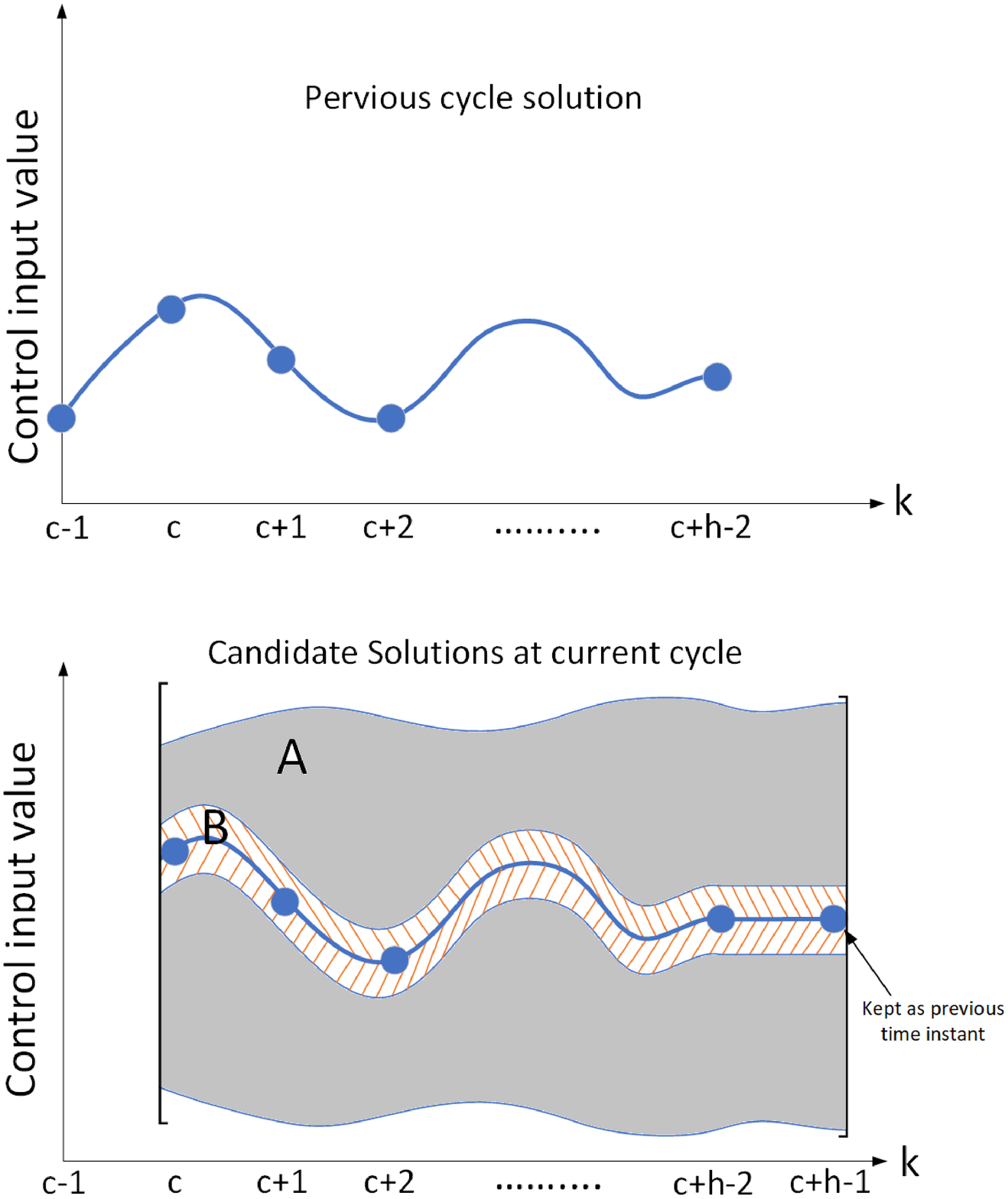}}
\end{subfigure} 
	\caption{A Graphical representation of the search spaces of GA. In (a), we show a contour (A) that represents the search space that comes from the physical constraints of the system, and another contour (B) marked with dashed red lines represents the BSM. In (b), we show how the candidate control inputs solutions at the current cycle are obtained by time-shifting the previous cycle control inputs by one-time step and searching around it in the BSM region marked by dashed red lines.}
	\label{fig:margins}
\end{figure*}

As we discussed, the GA explores the search space to estimate the best control inputs for the NMPC optimization at each control cycle. The GA uses the notion of \emph{margin} to represent the search space size. The margin is a constrained region likely to contain the best inputs. A small margin reduces the chance of finding the best inputs, while a larger one increases the chance but requires more computation.

One way to set the margin is by using the system's physical constraints (area A in Fig.~\ref{fig:margins}), i.e., the difference between the maximum and minimum allowable values of each input. However, this may result in a large search space. Another method is to define the margin as a neighborhood around the best inputs found in the previous cycle, shifted by a one-time step (area B in the figure). Tightening this margin reduces the search space but may lower the likelihood of finding optimal inputs.

This paper focuses on finding the best smallest margin (BSM), which improves the chance of finding optimal inputs with minimal computation. Throughout the rest of the paper, we denote the BSM of the $c^{\text{th}}$ control cycle as $\BSM$ as we define next.

\subsubsection*{Definition 1: Best smallest margin $\BSM$}

\begin{itemize}

\begin{figure}[t]
	\centerline{\includegraphics[width=.6\textwidth]{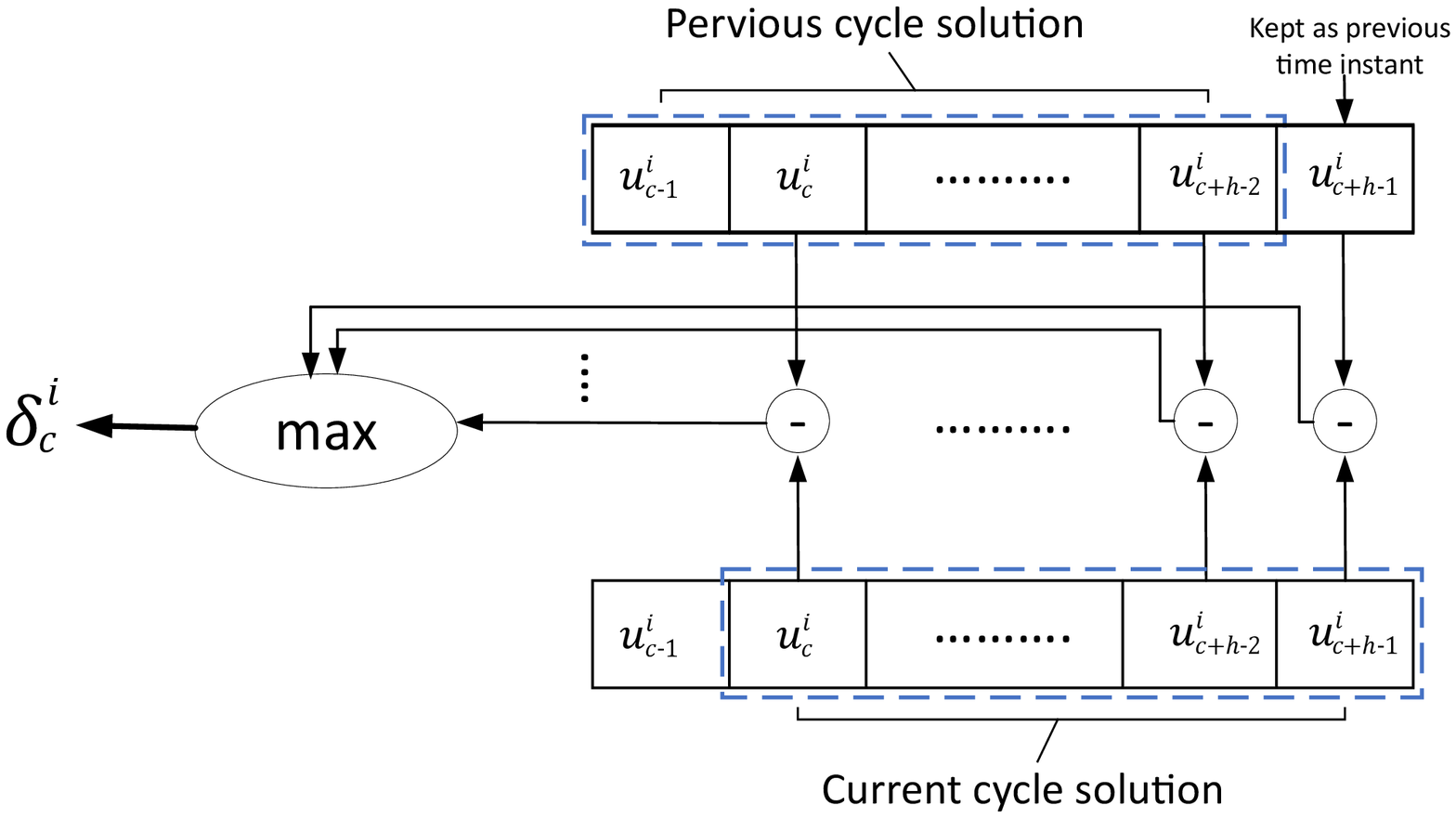}}
	\caption[Graphical illustration that presents the definition of $\bsm^i$.]{Graphical illustration that presents the definition of $\bsm^i$. We define $\bsm^i$ as the maximum of absolute differences between the corresponding elements of $\mathbf{g}^i_{c}$ that is estimated at the current $c^{\text{th}}$ and previous $c^{\text{th-1}}$ control cycles.}
	\label{fig:bsm}
\end{figure}

\item[] Let $\mathbf{g}^i_{c} = [u^i_{c}, \dots, u^i_{c+h-1}]$ be a vector that contains the values of the $i^{\text{th}}$ control input along the horizon $\hzn$ in the $c^{\textnormal{th}}$ control cycle. As depicted in Fig.~\ref{fig:bsm}, we define $\bsm^i$ as the maximum of absolute differences between the corresponding elements of $\mathbf{g}^i_{c}$ that is estimated at the current $c^{\text{th}}$ and previous $c^{\text{th-1}}$ control cycles. Then, the BSM vector is defined as 
\begin{equation}
\BSM  = [\bsm^1, \dots, \bsm^n].
\label{eq:BSM}
\end{equation}

\end{itemize}
In the following, we present the proposed approach to adaptively estimating the BSM. First, we will motivate the proposed approach by considering the main factors that affect it. Then, we present the proposed approach in detail.

\subsection{Factors affecting the BSM} 
\label{ssec:motivation}

In NMPC, the control inputs are estimated in each control cycle for a future horizon $h$. However, instead of time-shifting these readily computed control inputs to the next cycle, the NMPC estimates new control inputs to cope with the system's changing conditions. This is mainly because the system's future behavior may deviate from the expected behavior due to various factors, such as noise in both control inputs and state measurements. These noises result from inaccuracies in the actions applied to the system, model uncertainties, external disturbances, and changes in system dynamics.

Intuitively, suppose that noises are negligible. Then, in this case, the expected current system states $\mathbf{x}_c$ that is obtained by applying the control inputs estimated from the previous cycle in the mathematical model of the system~\eqref{eqn:statespaceform}, will slightly differ from the measured states in the current cycle $\mathbf{\bar{x}}_c$. Thus, we expect the control inputs at the current cycle to differ slightly from the previously estimated (time-shifted) inputs in the previous cycle. This means that the values of the BSM will be small in this case. In contrast, if these noises are considerable, then $\mathbf{x}_c$ will have significant differences from $\mathbf{\bar{x}}_c$. Thus, the system needs to aggressively change the control inputs from the previously estimated control inputs in the previous cycle. Therefore, larger values for the BSM are required.

To validate the above claim, we perform the following experiment. We use the UAV model discussed in Sec.I in the supplementary material and synthetically add errors in each control cycle to represent noise in control inputs and state measurements. The control input noise reflects inaccuracies in the actions applied to the system, while measurement noise captures errors arising from sensor inaccuracies, quantization, random fluctuations in sensor readings, and external disturbances. Then, we estimate the best control inputs in each cycle using the GA but with a large margin equal to the UAV model's physical constraints. Additionally, to neutralize the effect of the initial conditions, we started our experiments with initial state values equal to initial reference state values.

In each cycle, we obtain the associated BSM (as defined in Definition 1 and~\eqref{eq:BSM}) and record its maximum $\BSMmax = \max (\BSM)$. Also, we quantify the difference between $\mathbf{x}_c$ and $\mathbf{\bar{x}}_c$ by $\Ec$ defined as 
\begin{equation}
\Ec = [q^1(\Bar{x}^1_c-x^1_c)^2, \dots, q^m(\Bar{x}^m_c-x^m_c)^2],
\label{eq:error_vector}
\end{equation} 
and record its maximum $\Emax = \max (\Ec)$, where ${q}^{1},\dots,q^{m}$ weigh each element of the state vector according to their importance. We record the maximum values due to any deviation in any system state leading to a corresponding deviation in one or more control inputs because of the coupling between the inputs.

\begin{figure}[t]
	\centerline{\includegraphics[width=.6\textwidth]{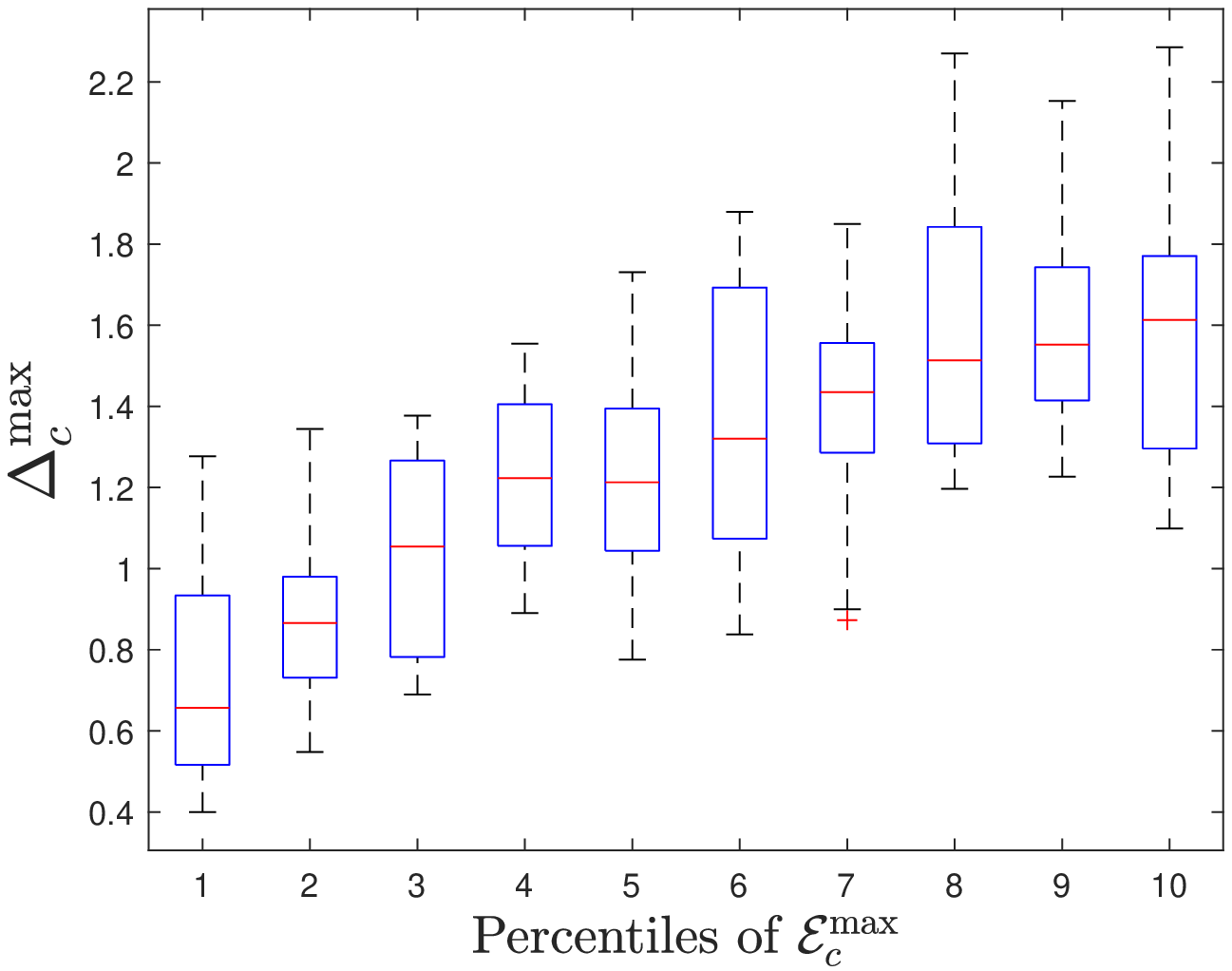}}
	\caption[A box plot shows the relation between the values of $\{\Emax\}$ and their corresponding $\BSMmax$ values.]{A box plot shows the relation between the values of $\{\Emax\}$ sorted in 10 percentiles and their corresponding $\BSMmax$ values. The $\BSMmax$ increases as $\Emax$ increases.}
	\label{fig:inpDiff}
\end{figure}

We ran the experiment several times with different reference trajectories and gathered the recorded $\Emax$ and $\BSMmax$ in all cycles. We sort the gathered values of $\Emax$ in 10 percentiles and plot the corresponding $\BSMmax$ against $\Emax$ using the box plots in Fig.~\ref{fig:inpDiff}. As shown in the figure, $\BSMmax$ increases with the increase in $\Emax$. Consequently, we can conclude that the $\BSMmax$ of a control cycle is proportionally related to the error between the expected and measured system states in the control cycle.

From the outcomes of the above experiment, one can estimate a mathematical relation that relates $\Emax$ with $\BSMmax$ and use this relation in the GA optimization to find the BSM to limit the search space at each cycle adaptively. However, the experiment is just for illustration and is far from realistic scenarios. This is because we used the maximum value among the elements of the BSM vector to report the results to neutralize the coupling effect between the control inputs. However, in real scenarios, we need to estimate the whole $\BSM$ vector in which each element is the corresponding control input's margin, as stated in Definition 1. 

In the next section, we present our proposed approach, which predicts $\BSM$ according to the error vector extracted in each cycle using a multivariate nonlinear support vector regression (MNSVR) algorithm.

\subsection{Proposed approach}
\label{sec:ProposedAGNMPC}

In this section, we introduce our proposed approach for estimating the best smallest margin ($\BSM$). The approach employs $n$ regression models to adaptively predict $\bsm^i$ for each control input during every control cycle. 

\begin{figure}[t]
	\centerline{\includegraphics[width=\textwidth]{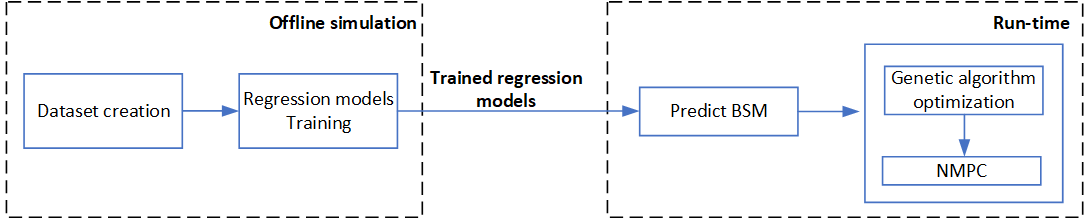}}
	\caption{Workflow of the proposed approach. The process begins with offline simulation, where a dataset is generated. Regression models are then trained on this dataset. In the runtime, the trained models predict the Best Smallest Margin (BSM), which is used for genetic algorithm optimization in NMPC.}
	\label{fig:propsed}
\end{figure}

As shown in Fig.~\ref{fig:propsed}, we start by building a synthetic dataset using the model of the system we want to control. The dataset creation procedure considers the inaccuracies in the actions applied to the system, model uncertainties, external disturbances, and changes in system dynamics. Then, we train the regression models on the dataset. Each regressor is trained to find the best parameters of the relationship\footnote{We drop the index $c$ from $\Ec$ and $\bsm$ for better readability.} between $\E$ and $\bsmm^i$. Specifically, the $i^{\text{th}}$ regressor finds the hyperplane $\mathbf{\omegaB}_i^T \phi(\E) + b_i$ that fits the maximum number of all $\bsmm^i$ in the dataset within a margin $\lambda$ from the hyperplane, where $\phi(\E)$ is a nonlinear transformation that maps $\E$ to a higher-dimensional space. In runtime, we use the trained regression models to estimate $\bsm^i$ for each control input and form the best smallest margin $\BSM$ from Eq.~\ref{eq:BSM}.

In the following, we give a detailed explanation of the regression models' training and prediction methodology. Next, we describe how the predicted $\BSM$ is integrated into the population generation process for the genetic algorithm (GA). Lastly, we explain the creation of the dataset used for training.

\subsubsection*{Training}

\begin{figure}[t]
\centerline{\includegraphics[width=.85\textwidth]{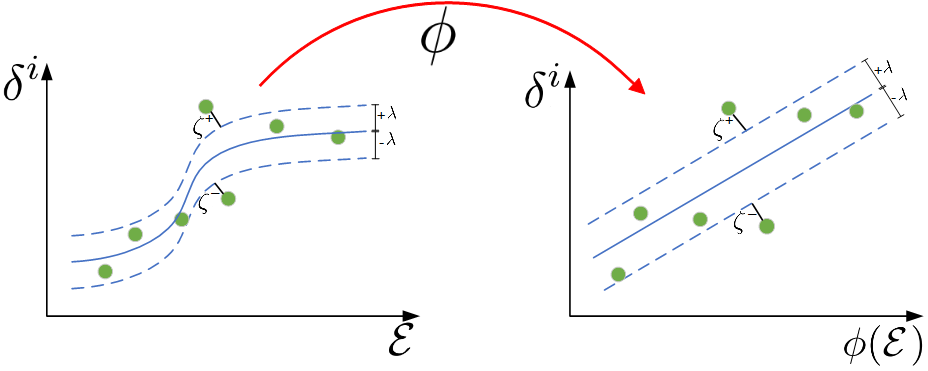}}
\caption{The function $\phi$ transforms the data into a higher dimensional feature space to make it possible to find a linear hyperplane that holds a maximum number of training observations within the margin.}
\label{fig:SVR}
\end{figure}

Given a training dataset comprising $D$ records, where the $j^{\text{th}}$ record in the dataset contains the error vector $\E_j$ and its associated $n$ values $\bsmm^1_j, \dots, \bsmm^n_j$, the goal of the training of the $i^{\text{th}}$ non-linear support vector regressor (SVR) model is to find the best parameters $\{\omegaB^*_{i}, b^{*}_{i} \}$ for the hyperplane by solving 
\begin{equation}
\begin{aligned}
  \{\omegaB^{*}_i, b^{*}_i \} =  &\arg \min_{\omegaB, b}   \frac{1}{2} \| \omegaB \| ^{2} + C \sum_{j=1}^{D} ( \zeta^{+}_j + \zeta^{-}_j ), \\
   &\text{s.t.}\\
   & \forall j:  \bsmm^i_j - \left(\omegaB^{T} \phi(\E_j)+b\right) \leq \lambda +  \zeta^{+}_j ,\\
   & \forall j:  \omegaB^{T}\phi(\E_j) + b - \bsmm^i_j \leq \lambda +  \zeta^{-}_j ,\\
   & \forall j: \zeta^{+}_j, \zeta^{-}_j \geq 0,
\end{aligned}
\label{eq:mnsvr}
\end{equation}
where $C > 0$ is a regularization parameter that penalizes the number of deviations larger than $\lambda$. The $\zeta^{+}_j$ and $\zeta^{-}_j$ are slack variables that allow the regression to have errors, as shown in Fig.~\ref{fig:SVR}. 

The dual optimization problem of \eqref{eq:mnsvr} is written as:
\begin{equation*}
\max_{\alpha_j^{+}, \alpha_j^{-}} \Bigg( -\frac{1}{2} \sum_{j=1}^{D} \sum_{k=1}^{D} (\alpha_j^{+} - \alpha_j^{-})(\alpha_k^{+} - \alpha_k^{-}) K(\E_j, \E_k) \quad - \lambda \sum_{j=1}^{D} (\alpha_j^{+} + \alpha_j^{-}) + \sum_{j=1}^{D} \bsmm^i_j (\alpha_j^{+} - \alpha_j^{-}) \Bigg),
\end{equation*}
s.t.
\begin{align*}
0 \leq \alpha_j^{+}, \alpha_j^{-} \leq C,\\
\sum_{j=1}^{D} (\alpha_j^{+} - \alpha_j^{-}) = 0,
\end{align*}
where \( \alpha_j^{+} , \alpha_j^{-} \) are Lagrange multipliers and \( K(\mathbf{x}, \mathbf{y}) = \phi(\mathbf{x})^{T}\phi(\mathbf{y})\) is the kernel function. We used a Gaussian kernel $K(\mathbf{x}, \mathbf{y}) = \exp\left( -\frac{1}{\kernelScale}\|\mathbf{x} - \mathbf{y} \|^2 \right)$ in this work, where $\kernelScale$ is the scaling factor which controls the shape of the decision boundary.

\subsubsection*{Prediction}
\begin{figure}[t]
	\centerline{\includegraphics[width=.6\textwidth]{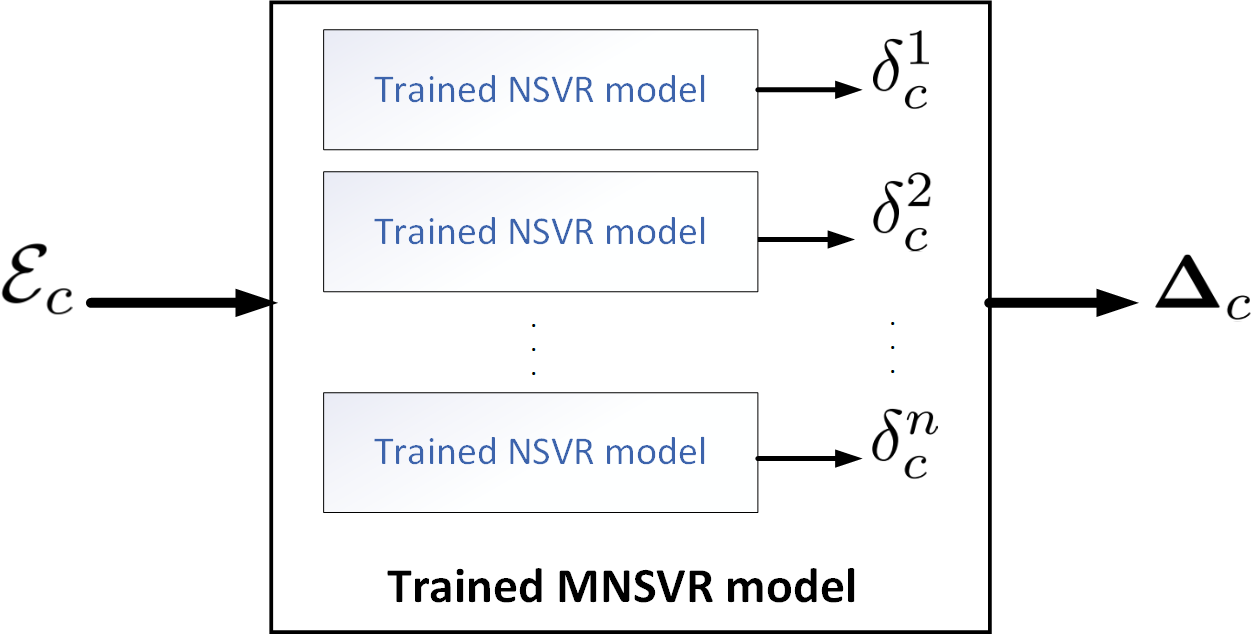}}
	\caption{Block diagram of the estimation of $\BSM$ in run time using the trained multivariate nonlinear support vector regressors (MNSVR).}
	\label{fig:SVM}
\end{figure}

Once the dual problem is solved in the training phase, the MNSVR is ready to make predictions. In the runtime, as depicted in Fig.~\ref{fig:SVM}, each regression model receives the input vector $\Ec$ for each control cycle and estimates
\begin{equation}
\bsm^i = \omegaB^{*T}_{i}\phi(\Ec) + b^{*}_{i}, \; \forall i=1, \dots, n.
\label{eq:runtime.svm}
\end{equation}
Or equivalently, by the dual form,
\begin{equation}
\bsm^i = \sum_{j\in \text{SV}} (\alpha_j^{+} - \alpha_j^{-}) K(\E_j, \Ec) + b^{*}_{i}, \; \forall i=1, \dots, n,
\label{eq:runtime.svm.dual}
\end{equation}
where $\text{SV}$ denotes the support vectors and \( b^{*}_{i} \) is computed as:
\[
b^{*}_{i} = \frac{1}{n_s} \sum_{j\in \text{SV}} \left( \bsmm^i_j - \sum_{j=1}^{D} (\alpha_j^{+} - \alpha_j^{-}) K(\E_j, \Ec) - \lambda\right),
\]
where $n_s$ denotes the number of support vectors. Please note that the dual form in~\eqref{eq:runtime.svm.dual} is preferred over \eqref{eq:runtime.svm} because of two reasons:
\begin{enumerate}
\item We do not explicitly compute the feature mapping $\phi(\E_c)$ because it may map to a very high (or even infinite) dimensional space. Instead, we use the kernel trick, which allows us to compute the inner product $\phi(\E_j)^{T}\phi(\E_c)$ directly using a kernel function without ever needing to calculate $\phi(\E_c)$.

\item The calculations in~\eqref{eq:runtime.svm.dual} are very fast because they only leverage the support vectors for making predictions. Specifically, support vectors are the only points with non-zero Lagrange multipliers $\alpha_j^{+}, \alpha_j^{-}$. Other points have their corresponding Lagrange multipliers equal to zero. Thus, the sum in~\eqref{eq:runtime.svm.dual} is calculated for a few points.
\end{enumerate}

In addition to predicting \(\bsm^i\) for all inputs, we estimate the confidence score of the prediction, $\mathcal{C}_c^i$, based on the distance between the input \(\Ec\) and the support vectors of the $i^{\text{th}}$ regressor. These confidence scores are then averaged to compute the overall prediction confidence, $\mathcal{O}_c$, which is used in the population generation process, as we detailed next. If the confidence is high, the prediction is used; otherwise, the margin is set to the physical boundaries due to unreliable estimation. Mathematically, we define
\begin{equation}
    \mathcal{C}_c^i = \frac{1}{\sum_{j\in \text{SV}} \|\boldsymbol{\mathcal{E}}_c - \mathbf{x}^{i}_j\|},
\end{equation}
where \(\mathbf{x}^{i}_j\) denotes the \(j^{\text{th}}\) support vector for the \(i^{\text{th}}\) regression model. The overall confidence score is then computed as
\begin{equation}
    \mathcal{O}_c = \frac{1}{n} \sum_{i=1}^{n} \mathcal{C}^i.
    \label{eq:overall_conf}
\end{equation}

\subsubsection*{Population generation}
Once we estimate $\bsm^i$ for all inputs, we can form the best smallest margin $\BSM$ according to~\eqref{eq:BSM}. Then, we generate a population of candidate solutions for the GA by sampling $p_c$ solutions randomly either within the estimated $\BSM$ or the physical margins of the control inputs, according to the quality of the estimated control inputs by the GA at the previous cycle. Specifically, if the previously estimated control inputs were of poor quality, the sampling is performed within the physical margins. This way, we ensure that the GA does not get confined to a search area around a sub-optimal solution from the previous cycle, preventing it from negatively impacting subsequent cycles' solutions.

Let $\finalBSM = [\psi_c^1, \dots, \psi_c^n]$ be the margin that the GA uses for sampling the $p_c$ solutions, which is defined as
\begin{equation}
\finalBSM = \begin{cases}
[\bsm^1, \dots, \bsm^n] & \text{if} \left( J(\mathbf{z}^{*}_{c-1}) \le J(\mathbf{z}^{*}_{c-2}) \,\textbf{or}\, J(\mathbf{z}^{*}_{c-1})<\epsilon \right) \,\textbf{and} \, \mathcal{O}_c> \eta,\\
[\beta^1,\dots,\beta^n] & \text{otherwise,}
\end{cases}
\label{eq:final-bsm}
\end{equation}
where $\beta^i$ is the physical margin of the $i^{\text{th}}$ control input (the difference between the maximum and minimum allowable values of this input), and $\epsilon$ and $\eta$ are predefined thresholds. According to~\eqref{eq:final-bsm}, the margin used by the GA for sampling the solutions is set to the estimated $\BSM$ when two conditions are met: (1) the estimated control inputs at the $c^{\text{th}-1}$ cycle are either better than those from the $c^{\text{th}-2}$ cycle or meet the cost function threshold $\epsilon$; and (2) the confidence of the margin prediction by our trained MNSVR exceeds the threshold $\eta$. Note that the thresholds $\epsilon$ and $\eta$ are set depending on the specific application and the desired level of performance.

In order to make $p_c$ independent on the scale of a specific control input, we normalize $\finalBSM$ according to the relative margin of each input. Specifically, we calculate the normalized margin $\bar{\bsm}^i$ for the $i^{\text{th}}$ control input as $\bar{\psi_c}^i = \psi_c^i / \beta^i$. Then, we compute the maximum of all normalized margins, denoted as $\alpha_c = \max \{ \bar{\psi_c}^i \}$. Finally, we set $p_c$ as
\begin{equation}
    p_c = \max\{\lfloor \SampDens  \alpha_c \rfloor,\xi\},
    \label{eq:pop-calc}
\end{equation}
where $\SampDens$ represents the maximum population size and $\xi$ represents a minimum population size, ensuring we generate enough candidate solutions even within small estimated margins. Once a population of $p_c$ candidate solutions is generated, GA proceeds normally to determine the best solution, as discussed in Sec.~\ref{sec:ga-npmc-solver}.

\subsubsection*{Dataset creation}
We followed the dataset creation method outlined in \cite{osa2018alDataset} to construct our training dataset, ensuring comprehensive coverage of various control scenarios. Specifically, we employed simulation techniques to systematically vary the system's control inputs within its constraints, generating numerous $N$ cycles of reference states. Consistent with Sec.\ref{ssec:motivation}, synthetic errors were introduced in each control cycle to represent noise in both control inputs and state measurements. Control input noise reflects inaccuracies in the actions applied to the system, while measurement noise captures errors arising from sensor inaccuracies, quantization, random fluctuations in sensor readings, and external disturbances.

We model the control input noise as an additive Gaussian noise $ \mathcal{N}(0, \Sigma_{\text{input}})$, where $\Sigma_{\text{input}} = \text{Diag}(\sigma^2_1,\dots,\sigma^2_n)$ is an $n\times n$ diagonal covariance matrix. Each $\sigma_i$ is defined as a fixed percentage $ \rho $ of the physical margin on the $i^{\text{th}}$ control input. Similarly, measurement noise is modeled as additive Gaussian noise $ \mathcal{N}(0, \Sigma_{\text{measure}})$ where $\Sigma_{\text{measure}}$ is an $m\times m$ diagonal covariance matrix with each diagonal element is set to a fixed percentage $\theta$ of the range of the corresponding state, i.e., the difference between the maximum and minimum allowable values. These synthetic errors ensure a more robust and realistic dataset for building our regression models.

In each cycle, control inputs are estimated using a genetic algorithm (GA) with $p_c$ solutions randomly sampled within the physical margins of the control inputs. A large $p_c$ and a high number of generations are used without a termination condition, allowing the GA to run until the end of its generations, thereby increasing the likelihood of converging to the best solution. This ensures that  $\bsm^i$ for each input is accurately estimated. As a result, each record in the dataset consists of the error vector $\Ec$ for the $c^\text{th}$ control cycle and its corresponding $n$ values of $\bsm^i$.

\subsection{Feasibility and stability}
\label{sec:feasibility_stability}
We ensure that the proposed approach maintains the feasibility of the control actions over the prediction horizon and leverages the inherent stability properties of the terminal state region $X_f$, as discussed in the assumptions A1-A4 in Sec.~\ref{sec:NMPC_formulation}.

Specifically, we follow a simple strategy to maintain the feasibility of the solution obtained with the proposed approach. We evaluate the feasibility of each candidate's solution in the population and exclude any infeasible ones. This ensures that the selected solution from the population afterward is feasible. Suppose the population has no feasible solution during any cycle or the computation time $t_c$ reaches a certain percentage of the maximum allowed control cycle's time. In that case, we retain the best-obtained solution from the previous cycle, shifted by a one-time step. 

Similarly, to ensure stability, we ensure that the terminal states at each control cycle lie within the designated terminal set $X_f$. By guaranteeing that the terminal state resides in this $X_f$ set, the system can be consistently guided toward the desired equilibrium. When the population has no solution that satisfies this condition or the computation time $t_c$ reaches a certain percentage of the maximum allowed control cycle's time, we retain the best-obtained solution from the previous cycle, shifted by a one-time step.

\subsection{Mathematical analysis and computational complexity}
\label{ssec:math}

In this section, we prove two key claims regarding the proposed approach. First, we prove that the probability of finding the best solution in one generation is higher in the search space encapsulated by the best smallest margin $\BSM$ provided by the proposed approach than in the larger search space encapsulating the system's physical constraints. Second, we show that the time required to find the best solution is greater in the larger search space encapsulating the system's physical constraints than in the smaller one encapsulated by the best smallest margin provided by the proposed approach. Additionally, we present the computational complexity of the proposed approach beyond the original genetic algorithm's complexity.

Let \( S \) represent the size of the original search space encapsulates the system's physical constraints, and \( S' \) the size of the reduced search space encapsulated by the best smallest margin provided by the proposed approach, where \( S' \leq S \). Let \( p_S \) and \( p_{S'} \) denote the population sizes for the spaces \( S \) and \( S' \), respectively. We assume that both search spaces, \( S \) and \( S' \), contain the best solution. Also, since \( S' \leq S \), it follows that \( p_{S'} \leq p_S \).

For each generation of solutions, the probability of \textbf{not} finding the best solution is:

\begin{itemize}
    \item For the search space of size \( S \): $\left( 1 - \frac{1}{S} \right)^{p_{S}}.$
    \item For the search space of size \( S' \): $\left( 1 - \frac{1}{S'} \right)^{p_{S'}}.$
\end{itemize}
Thus, the probability of finding the optimal solution in one generation is:
\begin{itemize}
    \item For search space \( S \): $ 1 - \left( 1 - \frac{1}{S} \right)^{p_{S}}. $
    \item For search space \( S' \): $ 1 - \left( 1 - \frac{1}{S'} \right)^{p_{S'}}. $
\end{itemize}
Since \( S' \leq S \), it follows that:
\[
\left( 1 - \frac{1}{S'} \right)^{p_{S'}} \leq \left( 1 - \frac{1}{S} \right)^{p_{S}}.
\]
This implies:
\begin{equation}
1 - \left( 1 - \frac{1}{S'} \right)^{p_{S'}} \geq 1 - \left( 1 - \frac{1}{S} \right)^{p_{S}}.
\label{eq:prob_proof}
\end{equation}
Thus, the probability of finding the optimal solution in one generation is higher in the smaller search space \( S' \) than in the larger search space \( S \). The inequality \eqref{eq:prob_proof} holds even if the population sizes for both spaces are forced to be the same.

Now, we show that the time needed to find the optimal solution is greater in the larger search space \( S \) than in the smaller search space \( S' \). There are two main reasons for this:
\begin{enumerate}
    \item Population size and exploration time: \\
    Since \( S' \leq S \), it follows that \( p_{S'} \leq p_S \). As a result, the time required to explore the population in the smaller search space \( S' \) is less than the time required in the larger search space \( S \).
    
    \item Higher probability of early termination: \\
    The inequality \eqref{eq:prob_proof} shows that the probability of finding the optimal solution is higher for the smaller search space \( S' \). Consequently, the probability of hitting a termination condition (which signals the end of the search) is also higher in the smaller search space, leading to faster termination compared to the larger search space.
\end{enumerate}

Finally, we present the computational complexity of the prediction phase of the multivariate nonlinear support vector regressor (MNSVR) used to predict the best smallest margin $\BSM$. This prediction phase is the additional overhead encompassed by the proposed approach compared to the original genetic algorithm. The prediction phase constitutes two operations: the prediction made using Eq.~\eqref{eq:runtime.svm.dual} and the confidence calculations with Eq.~\eqref{eq:overall_conf}. The calculations in Eq.~\eqref{eq:runtime.svm.dual} are mainly associated with kernel computation for each support vector. For the Gaussian kernel used in our work, it involves computing the Euclidean distance between $ \E_c $ and each support vector, which requires $ O(m) $ operations per support vector (because the size of $\E_c$ is $m$). Then, for $n_s$ support vectors, the complexity of prediction per control input is $O(m \cdot n_s)$. Similarly, the confidence calculations in Eq.~\eqref{eq:overall_conf} involve computing the Euclidean distance between $ \E_c $ and each support vector. So, it has a complexity of $O(m \cdot n_s)$. Thus, the overall prediction complexity for $n$ inputs with the MNSVR is:
\begin{equation}
O(2n \cdot m \cdot n_s)
\nonumber\end{equation}

The complete algorithm of the proposed approach is outlined in Algorithm~\ref{alg:wholeprocess}.

\begin{algorithm}[ht!]
\begin{algorithmic}[1]
\item[] $N$: is the number of generated reference states for training the MNSVR.
\item[] $H$: is the number of control cycles.
\item[] $G$: is the number of generations.
\item[] $\upsilon$: is the percentage of the maximum allowed control cycle's time.
	\State Generate $N$ cycles of reference states.
	\State $c$ $\gets$ 0.
    \While{$c<N$}
    \State Compute $\Ec$ using~\eqref{eq:error_vector}, and its associated $n$ values of $\bsm^i$.
    \State  Add $\Ec$ along with its associated $n$ values of $\bsm^i$ to the dataset records.
    \State $c \gets c+1$.
    \EndWhile
	\State Train $n$ MNSVR models to estimate $\{\omegaB^*_{i}, b^{*}_{i} \}$,  $\forall i=1, \dots, n$.
	\State $c$ $\gets$ 0.
	\While{$c< H$}
      \State Compute $\Ec$ using~\eqref{eq:error_vector}.
      \State Estimate $\bsm^i$ using~\eqref{eq:runtime.svm.dual}.
      \State Compute $\finalBSM$ using~\eqref{eq:final-bsm}.
      \State Compute $p_c$ using~\eqref{eq:pop-calc}.
      \State Sampling $p_c$ solutions from $\finalBSM$ search space.
      \State $g \gets 0$.
 	  \While{$(g<G)  \text{ and } (t_{c} \leq \upsilon t_{s})$}
              \State Calculate solutions fitness using~\eqref{eq:fitness}.
              \State Exclude infeasible solutions and solutions that lead to instability according to Sec.~\ref{sec:feasibility_stability}
 	     \State Sort solutions according to their fitness values.
 	     \State $\mathbf{z}_c^* \gets$ solution with the minimum fitness.
 	     \If {$J(\mathbf{z}_c^*) \le \epsilon$}
 	     \State \textbf{break.}
 	     \EndIf
 	     \State Apply selection, cross-over, and mutation operations.
 	     \State $g \gets g+1$.
	 \EndWhile
	 \State Retain the first $n$ elements from $\mathbf{z}_c^*$ as best solution $\mathbf{u}_c \gets \mathbf{z}_c^* [1:n]$.
      \If {$\mathbf{u}_c =\emptyset $}
         \State $\mathbf{u}_c \gets \mathbf{z}_{(c-1)}^* [n+1:2n]$
      \EndIf
	 \State Apply $\mathbf{u}_c$ to the plant and get the new states $\mathbf{\bar{x}}_{c+1}$.
	  \State $c \gets c+1$.
  	\EndWhile
\end{algorithmic}
\caption{Proposed learning of optimal search space size for genetic optimization.}
\label{alg:wholeprocess}
\end{algorithm}

\section{Experimental results}
\label{sec:results}

In this section, we experimentally assess the performance of the proposed approach in comparison to four evolutionary optimization techniques: Particle Swarm Optimization (PSO)~\cite{diwan2023PSOcompare}, Differential Evolution (DE)~\cite{zhang2021DEcompare}, Original Genetic (OG), and Modified Genetic (MG) algorithms~\cite{arrigoni2022trjplanningGAMPC,du2016develGAsaftyEqn}. Both OG and MG maintain a constant population size across control cycles. However, MG divides its population into two parts, one generated within a fixed physical margin (similar to OG) and the other comprising time-shifted versions of the best solutions from the previous control cycle's genetic optimization.

We conducted experiments using the proposed approach and compared methods to control three nonlinear systems: an Unmanned Aerial Vehicle (UAV), a ground vehicle, and a Single Fixed Joint Robot (SFJR). Sections I, II, and III of the supplementary material provide detailed descriptions of these systems. For all experiments, we assume the models' states are fully measurable and available to the NMPC.

To train the regression models used in our approach, we generated a dataset with $ D = 10,000,000 $ records for both the UAV and ground vehicle models and $ D = 800,000 $ for the SFJR model. The parameters of the regression models were estimated using a 5-fold cross-validation technique~\cite{chorowski2014crossvalidation}, with Mean Squared Error (MSE) and the correlation coefficient (R-value) used as performance metrics in the training phase~\cite{tatachar2021metrics}. The values of the hyper-parameters of the proposed approach are provided in Section S.IV of the supplementary material. In all approaches, relatively long horizons were used to ensure system stability, as recommended in~\cite{boccia2014longH_1,grimm2005longH_2,grune2010longH_3}.

Our comparison evaluates both performance and computational time. The performance metric is defined as the average cost across all control cycles. Let $H$ represent the number of control cycles (consistent across all experiments). The average cost, $E$, is calculated as:
\begin{equation}
    E = \frac{1}{H} \sum_{c=1}^{H} J(\mathbf{z}_c^*).
    \label{eq:avgCost}
\end{equation}

All approaches are implemented using MATLAB on the GPU of an Nvidia\textsuperscript{\texttrademark} Jetson TX2 embedded platform in a processor-in-the-loop (PIL) fashion~\cite{PILNvidiaURL}. The GPU implementations of GA, PSO, and DE utilized MATLAB's built-in GPU parallelization features. To ensure a fair comparison, the NMPC settings for each nonlinear system remained consistent across all optimization techniques during the experiments.

In the following, we first introduce the Processor-In-the-Loop (PIL) setup used in all our experiments, providing complete details of the hardware and software configurations. Next, we present the experiments that were conducted to evaluate and compare our proposed approach against the other evolutionary optimization algorithms. Afterward, we discuss the robustness of the proposed approach, particularly focusing on its generalizability and robustness against various uncertainties in real-world conditions supported by additional validation experiments. Finally, we present a parameter sensitivity analysis to evaluate the impact of key thresholds and parameters on the performance of our approach.

\subsection{Processor-in-the-loop setup}

\begin{figure}[t]
	\centerline{\includegraphics[width=.9\textwidth]{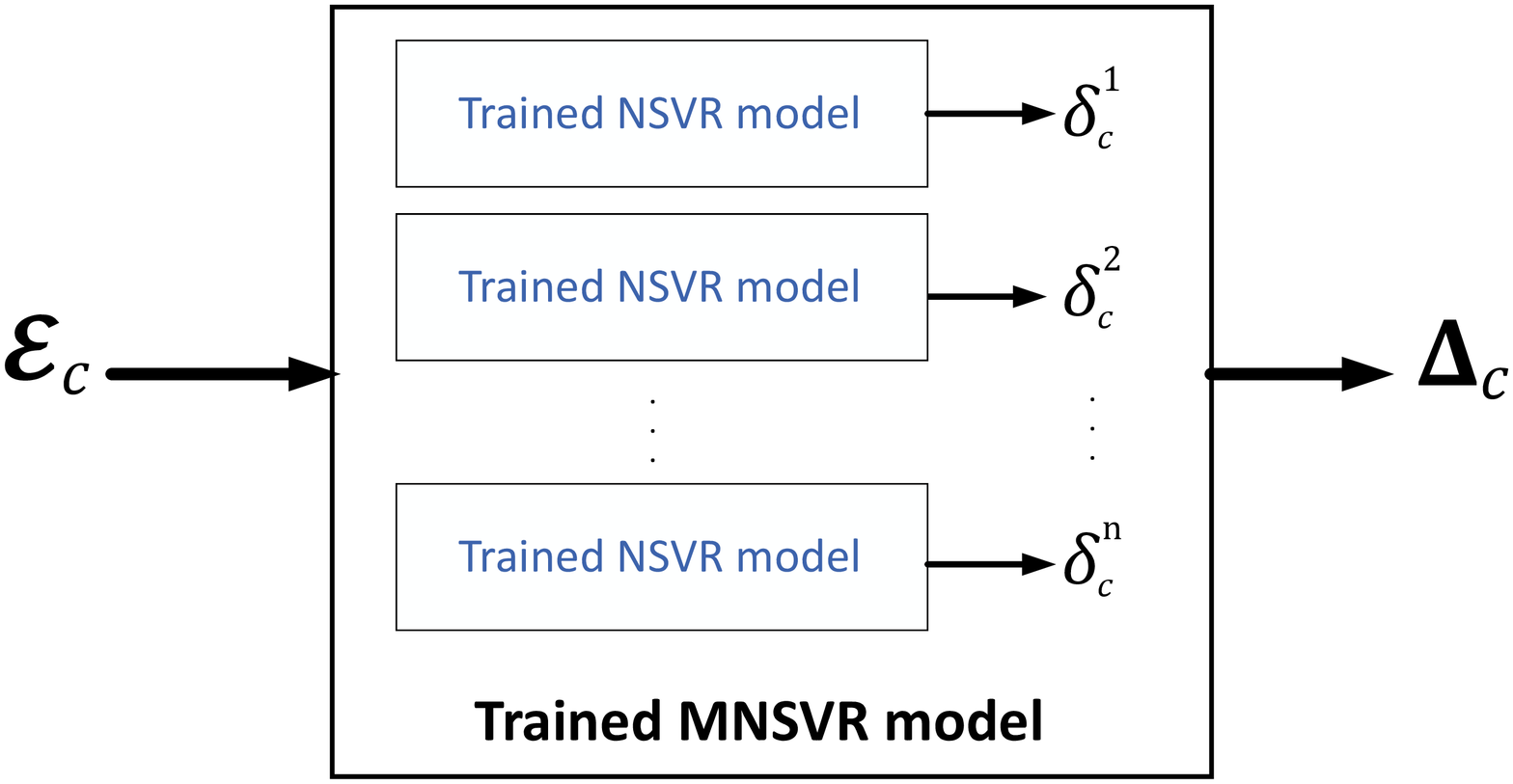}}
	\caption[The processor-in-the-loop setup.]{The processor-in-the-loop setup.}
	\label{fig:PIL}
\end{figure}

As illustrated in Fig.~\ref{fig:PIL}, our hardware-in-the-loop setup consists of two components: a plant simulator running on a host machine and the NMPC controller running on an embedded target. These two components are connected via an RS232 Ethernet cable. The host machine sends reference inputs and states to the controller and measures the plant dynamics after applying the control inputs. The host machine is a laptop running Microsoft Windows 10 OS with an Intel-i7 8700 CPU and 16 GB of RAM. All implementations are written in MATLAB. The embedded target is an NVIDIA Jetson TX2, running Ubuntu 16.04.

For the UAV model, we aim to control its motion in 3D space. The UAV's plant simulator receives rotor velocities from the controller as control inputs and outputs position, angular velocities, linear velocities, and orientation data. For the ground vehicle, the focus is on controlling its longitudinal and lateral positions and orientation. The plant simulator accepts steering angle and acceleration as inputs and outputs longitudinal/lateral positions, velocities, yaw angle, and yaw rate. For the SFJR, the goal is to control the link angle using the DC motor voltage. The plant simulator provides angular positions, velocities (motor and link sides), and motor current as outputs.

In all experiments, reference states were designed with rapid and slow variations to simulate real-world conditions. The supplementary material illustrates these reference states for all models in Fig. S.3. The optimization hyper-parameters used for the UAV, ground vehicle, and SFJR models are detailed in Tables S.II, S.IV, and S.VI.

\subsection{Comparison with other approaches}
\label{sec:results-Exp1}
This section presents our experimental comparison between the proposed approach and the other four evolutionary algorithms: PSO, DE, OG, and MG. To ensure fairness in comparison for OG and MG approaches, we set the maximum population size $\SampDens$ in~\eqref{eq:pop-calc} to be the same value for both. Also, for OG and MG approaches, we set the margin to be the physical margin of each control input, i.e., we set $\bsm^i = \beta^i$ for $i=1, \dots, n$. This means that $\alpha_c$ in~\eqref{eq:pop-calc} has the value 1 for both. In the MG approach, 80\% of the population is randomly generated within the physical margin, while the remaining 20\% is obtained by shifting the best solutions from the previous control cycle backward by one-time step~\cite{arrigoni2022trjplanningGAMPC,du2016develGAsaftyEqn}.

We compute the average cost $E$ and the average computational time for all approaches when applying to control the UAV, the ground vehicle, and the SFJR models using the NMPC. Additionally, we compute the convergence rate, which is the percentage of the cycles that each approach converges to the optimal solution\footnote{In all our experiments optimality means $J(\mathbf{z}_c)<\epsilon$} before the termination of the cycle's time. Tabel~\ref{tab:EX1} presents the compared approaches' average computational time (in milliseconds), convergence rate, and average cost $E$ for each approach when applied to control the UAV, ground vehicle, and SFJR models using NMPC.

\begin{table}[t]
    \footnotesize
    \centering
    \caption{Comparison of the average computation time, convergence rate, and average cost $E$ for each approach when applied to control the UAV, ground vehicle, and SFJR models using NMPC. The best results are highlighted in \textbf{bold}, and the second-best results are marked in \textcolor{red}{red}. The last row shows the percentage improvement of the proposed approach compared to the second-best approach.}
    \renewcommand{\arraystretch}{1.5}
    \begin{tabularx}{\textwidth}{|c|*{3}{>{\centering\arraybackslash}X}|*{3}{>{\centering\arraybackslash}X}|*{3}{>{\centering\arraybackslash}X}|}
        \hline
        \multirow{2}{*}{Approach} & \multicolumn{3}{c|}{AVG. COMP [msec]} & \multicolumn{3}{c|}{Convergence [\%]} & \multicolumn{3}{c|}{E} \\ 
        \cline{2-10}
        & UAV & Vehicle & SFJR & UAV & Vehicle & SFJR & UAV & Vehicle & SFJR \\ 
        \hline
        PSO  & \textcolor{red}{16.6} & \textcolor{red}{27.5} & 38.8 & \textcolor{red}{40} & 33 & 28 & 1.4 & \textcolor{red}{4.2} & 3.9 \\
        DE   & 19.2 & 37.2 & \textcolor{red}{29.5} & 22 & 30 & \textcolor{red}{40} & 2.2 & 6.7 & 2.4 \\
        OG   & 18.8 & 38.9 & 38.9 & 26 & 30 & 30 & 1.5 & 5.7 & 2.7 \\
        MG  & 18.1 & 35.2 & 35.2 & 24 & \textcolor{red}{36} & 36 & \textcolor{red}{1.3} & 5.9 & \textcolor{red}{2.2} \\
        \textbf{Proposed}   & \textbf{13.8} & \textbf{19.2} & \textbf{16.2} & \textbf{62} & \textbf{66} & \textbf{75} & \textbf{0.7} & \textbf{3.4} & \textbf{1.4} \\
        \hline
        Improvements [\%] & 16.87 & 30.18 & 45.08 & 35.48 & 45.45 & 46.67 & 46.15 & 19.05 & 36.36 \\
        \hline
    \end{tabularx}
        \label{tab:EX1}
\end{table}
As shown in Table.~\ref{tab:EX1}, the proposed approach significantly reduces the average computational time compared to other approaches. Specifically, the reduction in the computation time compared to the second-best approach reaches 16\%, 30\%, and 45\% for the UAV, the vehicle, and the SFJR, respectively. This reduction in computational time results from exploring less search space using the proposed approach. Therefore, it converges to the optimal solution in a shorter time. This result is aligned with the mathematical analysis presented in Sec.~\ref{ssec:math}. Additionally, the proposed approach outperforms the other approaches in terms of the convergence rate. Specifically, the proposed approach shows superior convergence rates. Compared to the second-best approach, the proposed approach has a better convergence rate by 35\%, 45\%, and 46\% for the UAV, the vehicle, and the SFJR models, respectively. Again, the result is aligned with the mathematical analysis presented in Sec.~\ref{ssec:math}. Based on the inequality \eqref{eq:prob_proof}, the probability of finding the optimal solution is higher for the smaller search space provided by the proposed approach. Finally, the proposed approach outperforms other approaches in $E$. The table shows that the proposed approach reduces $E$ compared to the second-best approach by 46\%, 19\%, and 36\% for the UAV, vehicle, and SFJR, respectively. The reason is that the proposed approach explores less search space due to its estimation of the optimal margin, and thus, it converges to the optimal solution more often.

\begin{figure}[t]
	\centering
\begin{subfigure}[]{0.32\textwidth}
   \centerline{\includegraphics[width=\textwidth]{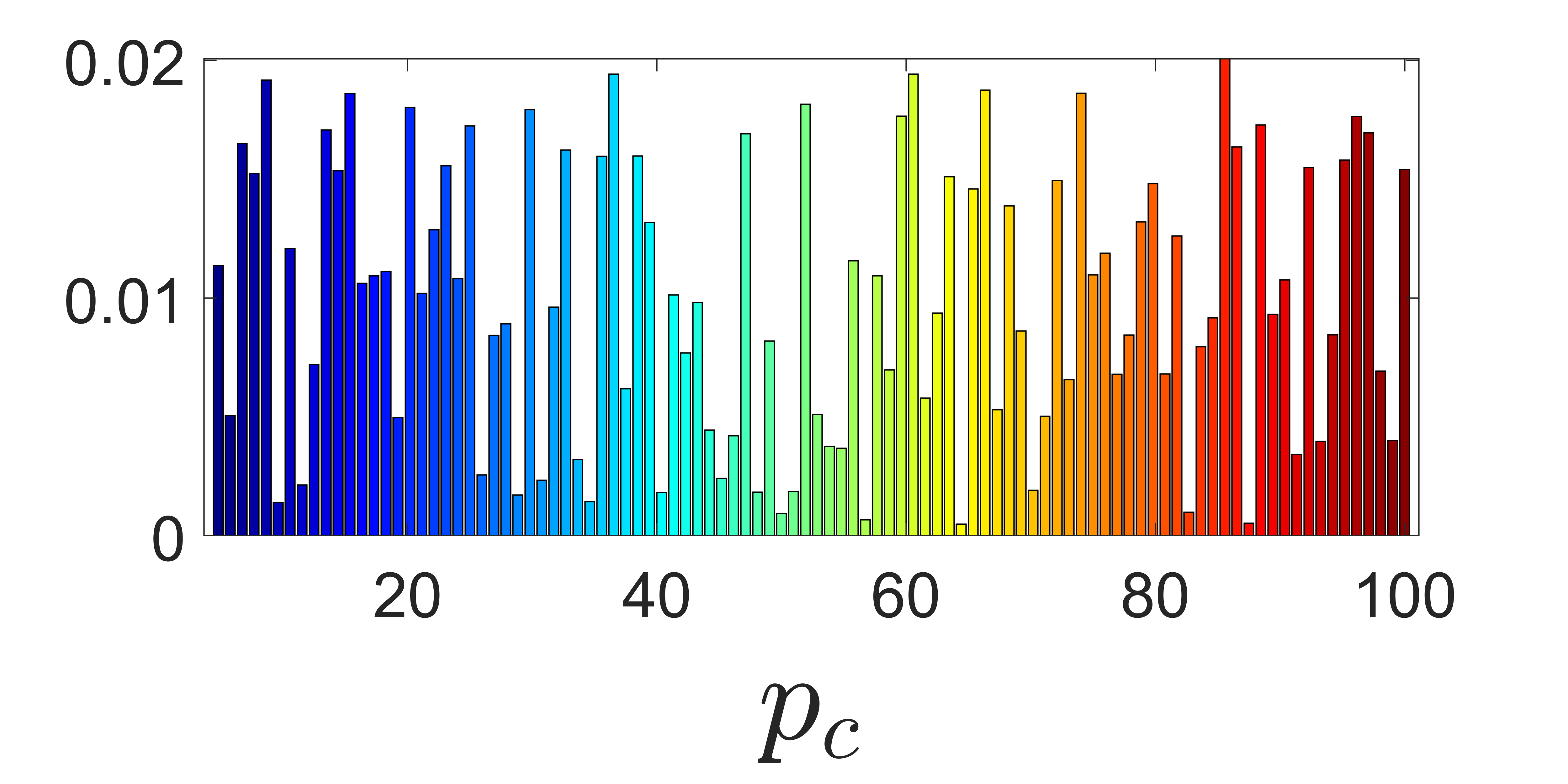}}
   \caption{}
\end{subfigure}
\begin{subfigure}[]{0.32\textwidth}
   \centerline{\includegraphics[width=\textwidth]{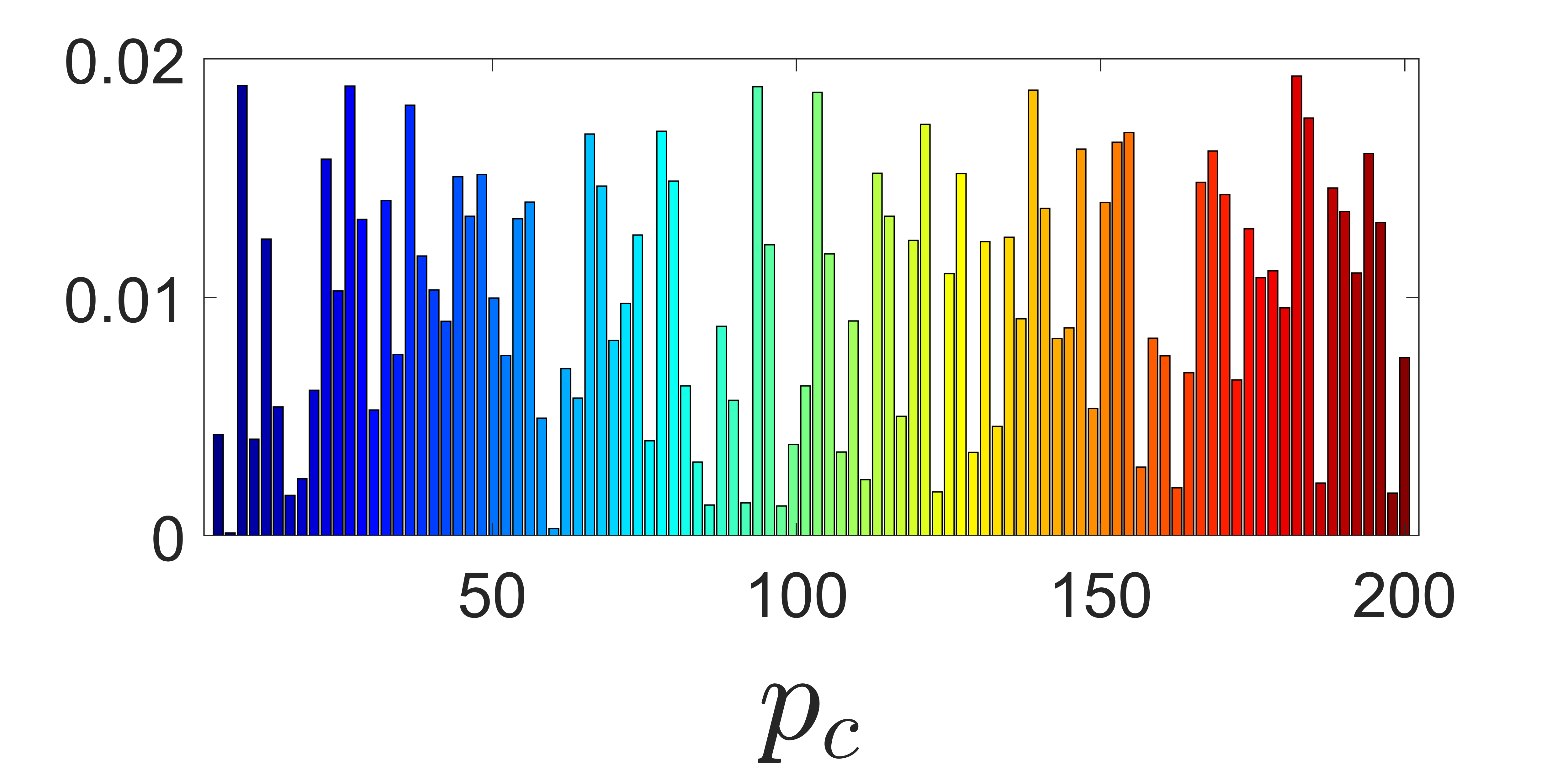}}
    \caption{}
\end{subfigure}
\begin{subfigure}[]{0.32\textwidth}
   \centerline{\includegraphics[width=\textwidth]{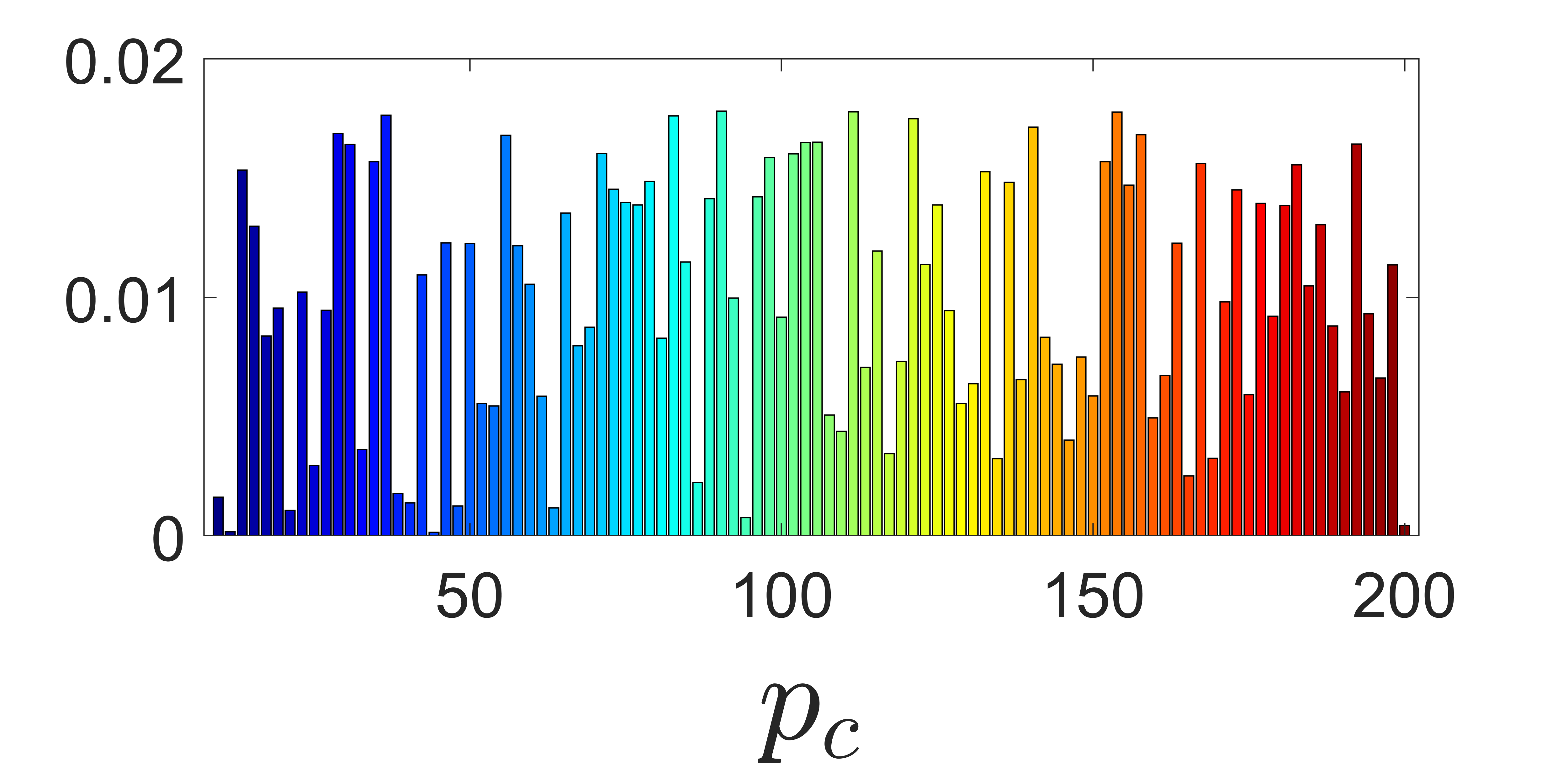}}
    \caption{}
\end{subfigure}
	\caption{Histogram of the population size $p_c$ in all cycles for the proposed approach. The proposed approach utilizes a variable population size according to the estimated margin instead of employing a fixed, large population size in all control cycles, which significantly reduces computations. Figures (a), (b), and (c) correspond to the UAV, ground vehicle, and SFJR models, respectively.}
	\label{fig:histogram}
\end{figure}

To further show the effectiveness of the proposed approach in reducing the population size $p_c$ in the conducted experiments, we plot a histogram for $p_c$ in all cycles for the UAV, the vehicle, and the SFJR models in Fig.~\ref{fig:histogram} (a), (b), and (c), respectively. The figure shows that the proposed approach exhibits variable $p_c$, which results from its adaptive estimation for the optimal search space. Additionally, $p_c$ has small values most of the time, as shown by the histograms, which significantly reduces computations. On the other hand, the OG and MG approaches utilize a fixed, large population to cover the large search space that arises from the physical margins of the control inputs, which results in unnecessary computations.

\subsection{Robustness}

The proposed approach relies on a synthetic dataset to train the prediction model to determine the best smallest margin. This dataset is generated in a simulated environment that may differ from real-world conditions, potentially introducing challenges in practical applications. In this section, we address the robustness of the proposed method and evaluate its performance in real-world scenarios.

As described in Sec.~\ref{sec:ProposedAGNMPC}, synthetic errors were introduced in each control cycle during dataset generation to represent noise in both control inputs and state measurements. These noises are characterized by their diagonal covariance matrices, $\Sigma_{\text{input}}$ and $\Sigma_{\text{measure}}$. The parameters $\rho$ and $\theta$ influence the entries of these matrices.

\begin{figure*}[t] 
\centering 
\begin{minipage}[]{.3\textwidth} 
\centerline{\includegraphics[width=\textwidth]{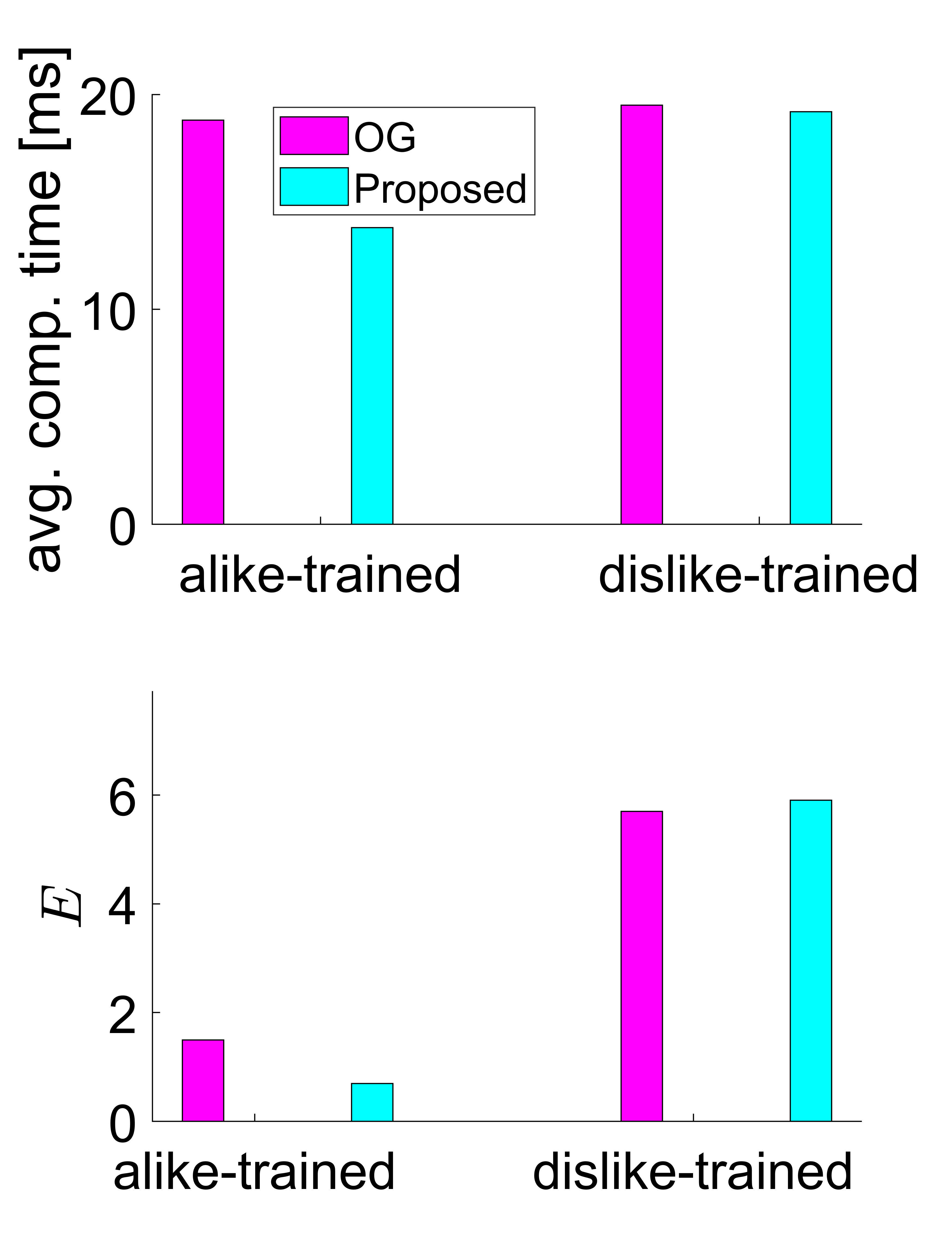}} 
           \captionsetup{labelformat=empty}
\caption{(a)} \end{minipage} 
\begin{minipage}[]{.3\textwidth} 
\centerline{\includegraphics[width=\textwidth]{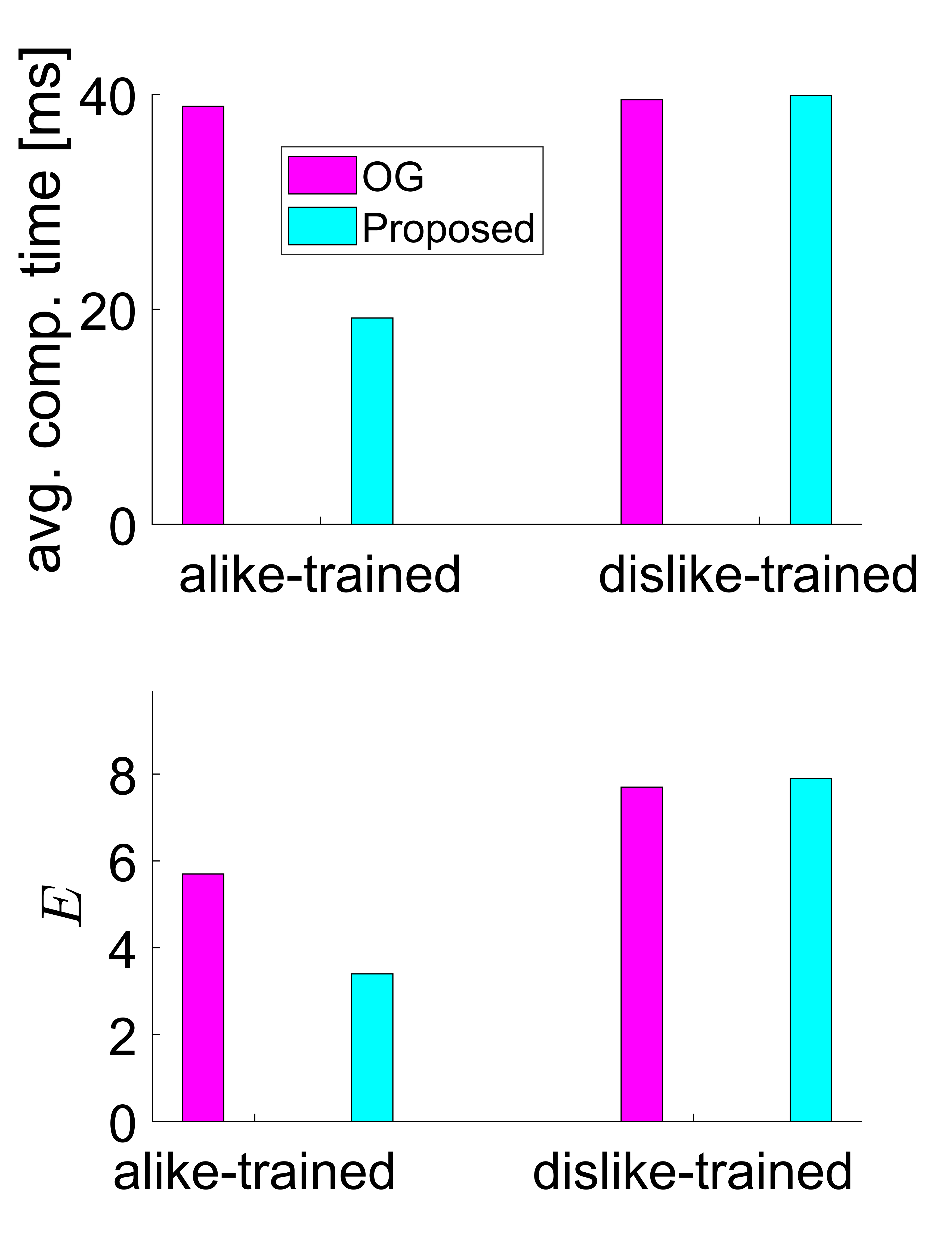}}
           \captionsetup{labelformat=empty}
\caption{(b)} 
\end{minipage}
\hspace*{0.1in}
\begin{minipage}[]{.3\textwidth} 
\centerline{\includegraphics[width=\textwidth]{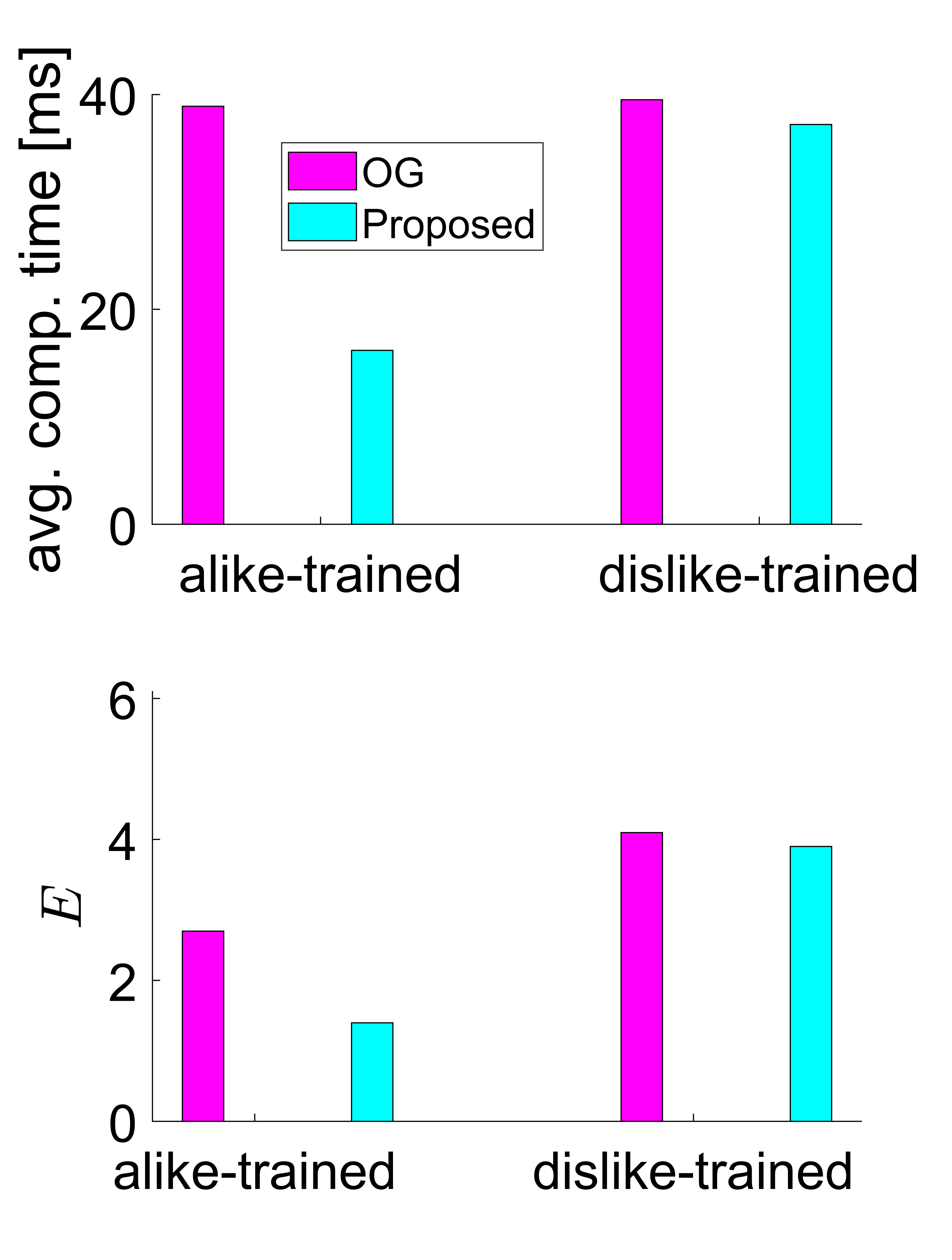}} 
           \captionsetup{labelformat=empty}
\caption{(c)}
\end{minipage} 
\caption{Performance comparison between the proposed approach and the OG in \textit{alike-trained} and \textit{dislike-trained} scenarios. The proposed approach shows a significant reduction in both computational time and average cost in the \textit{alike-trained} environment, leveraging effective predicted BSM for faster convergence. In the \textit{dislike-trained} scenario, the proposed approach maintains performance similar to OG by expanding the search space to the physical margins and preventing any degradation in performance relative to OG. Figures (a), (b), and (c) correspond to the UAV, ground vehicle, and SFJR models, respectively.}
\label{fig:FGcomparison} 
\end{figure*}

We performed the following experiment to evaluate the robustness of the proposed approach and demonstrate its performance in real-world scenarios. We evaluated the average cost ($E$) and average computational time of the proposed approach in two cases: (a) when $\rho$ and $\theta$ values differ from those used during dataset creation (referred to as \textit{dislike-trained}), and (b) when the values match those used during training (referred to as \textit{alike-trained}). In both cases, we treat the Original Genetic Algorithm (OG) as a baseline for comparison. We plot $E$ and average computational time in Fig.~\ref{fig:FGcomparison}.

As depicted in Fig.~\ref{fig:FGcomparison}, the proposed approach consistently outperforms the OG algorithm in the \textit{alike-trained} scenario, significantly reducing both computational time and average cost. This improvement stems from the proposed method's ability to effectively leverage the predicted BSM across many control cycles, allowing for faster convergence to optimal solutions while minimizing computational effort.

In the \textit{dislike-trained} scenario, however, the proposed approach performs similarly to OG. This occurs because the approach produces low-confidence predictions for the BSM in many cycles when dealing with noise in control inputs or state measurements outside its training range. In response, the proposed method adapts by expanding the genetic algorithm's search space to the physical margins, as described in~\eqref{eq:final-bsm}. This adaptive mechanism prevents reliance on inaccurate search space predictions, thereby avoiding significant performance degradation compared to OG. In other words, while the proposed approach outperforms OG in the \textit{alike-trained} scenario, its performance does not drop below OG's in the worst-case scenario when the real-world environment differs from the training conditions, i.e., in the \textit{dislike-trained} scenario.

\subsection{Parameter sensitivity analysis}
\begin{figure*}[t]
	\centering
\begin{subfigure}[]{0.31\textwidth}
   \centerline{\includegraphics[width=\textwidth]{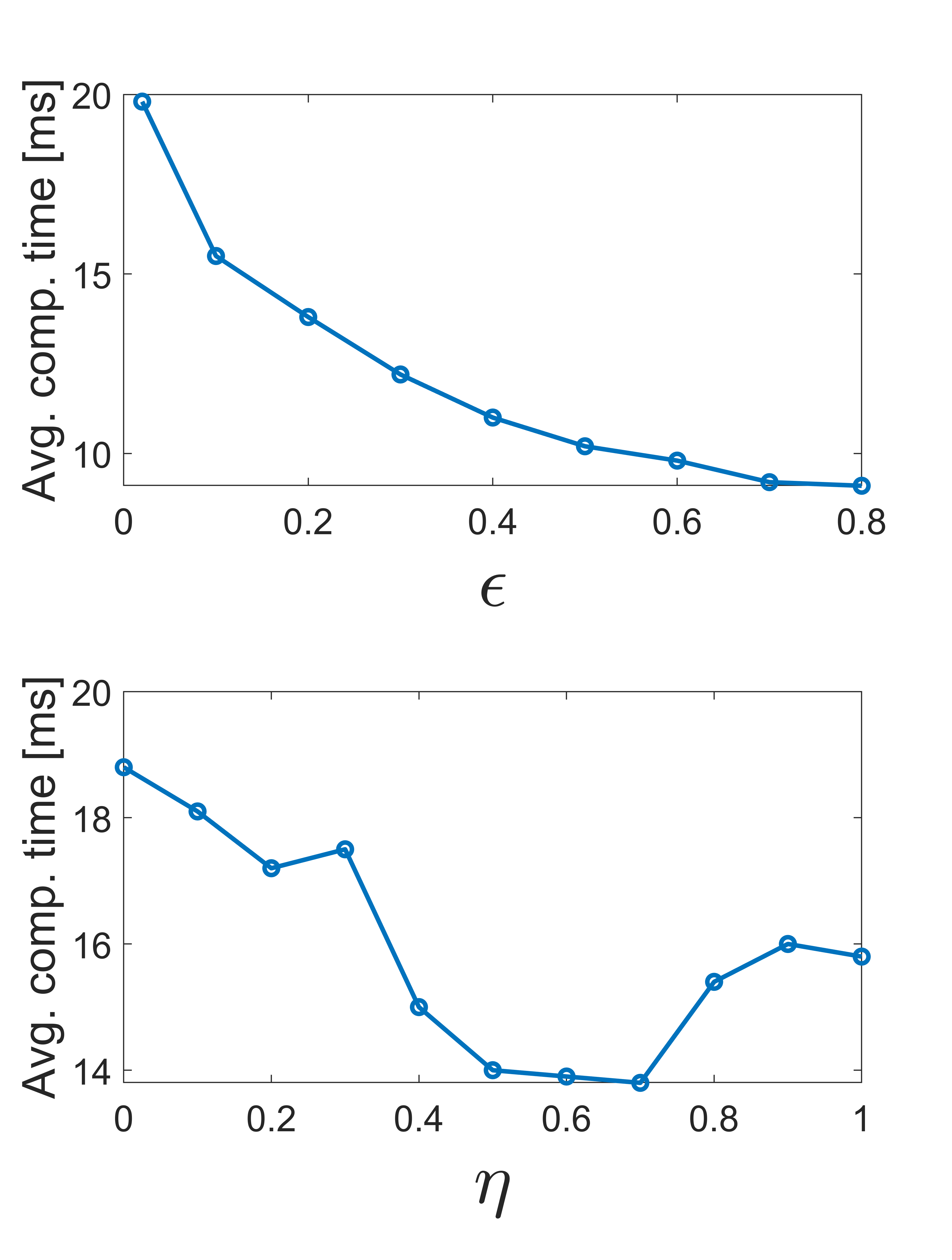}}
\end{subfigure}
\hspace*{0.1in}
\begin{subfigure}[]{0.31\textwidth}
   \centerline{\includegraphics[width=\textwidth]{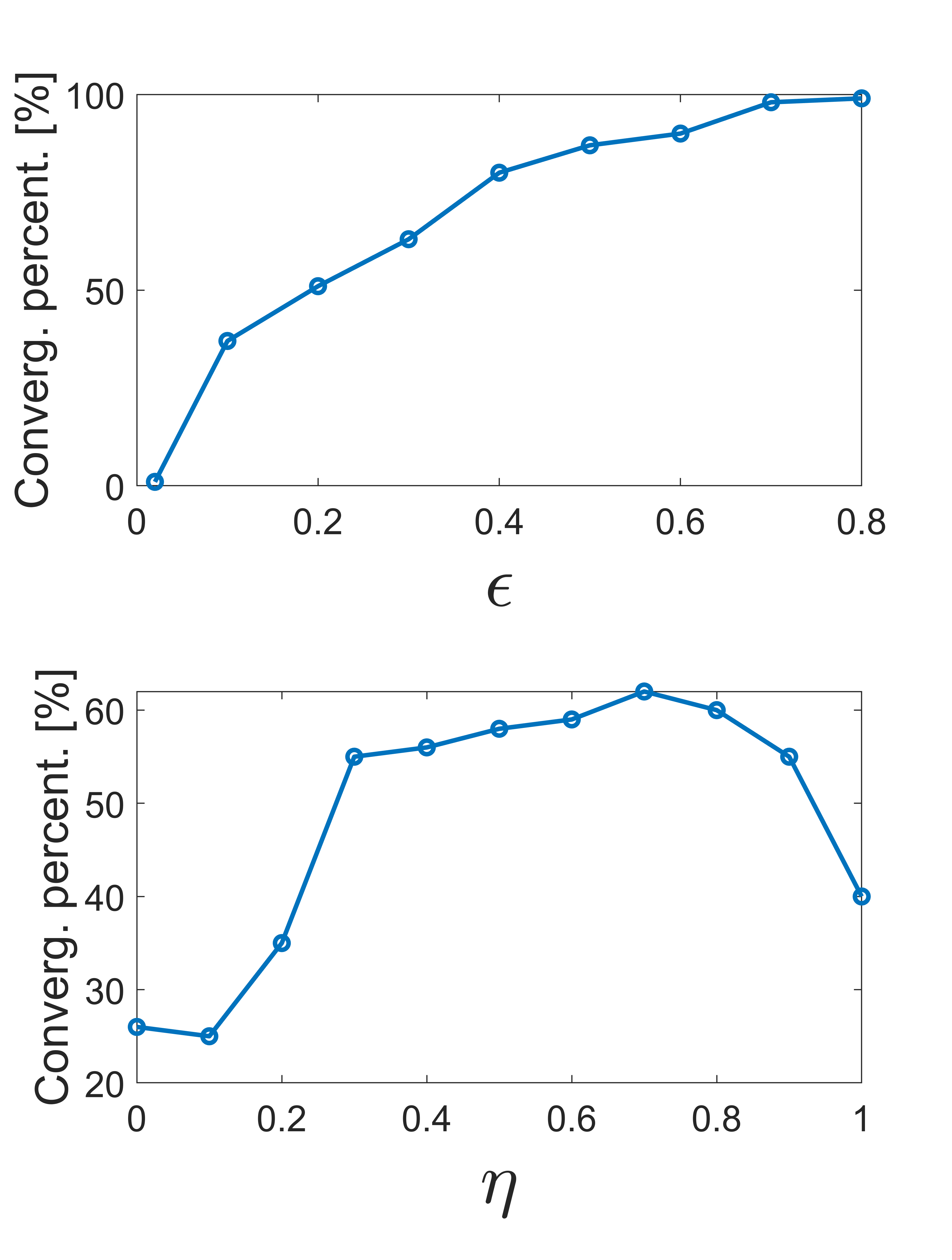}}
\end{subfigure}
\begin{subfigure}[]{0.31\textwidth}
   \centerline{\includegraphics[width=\textwidth]{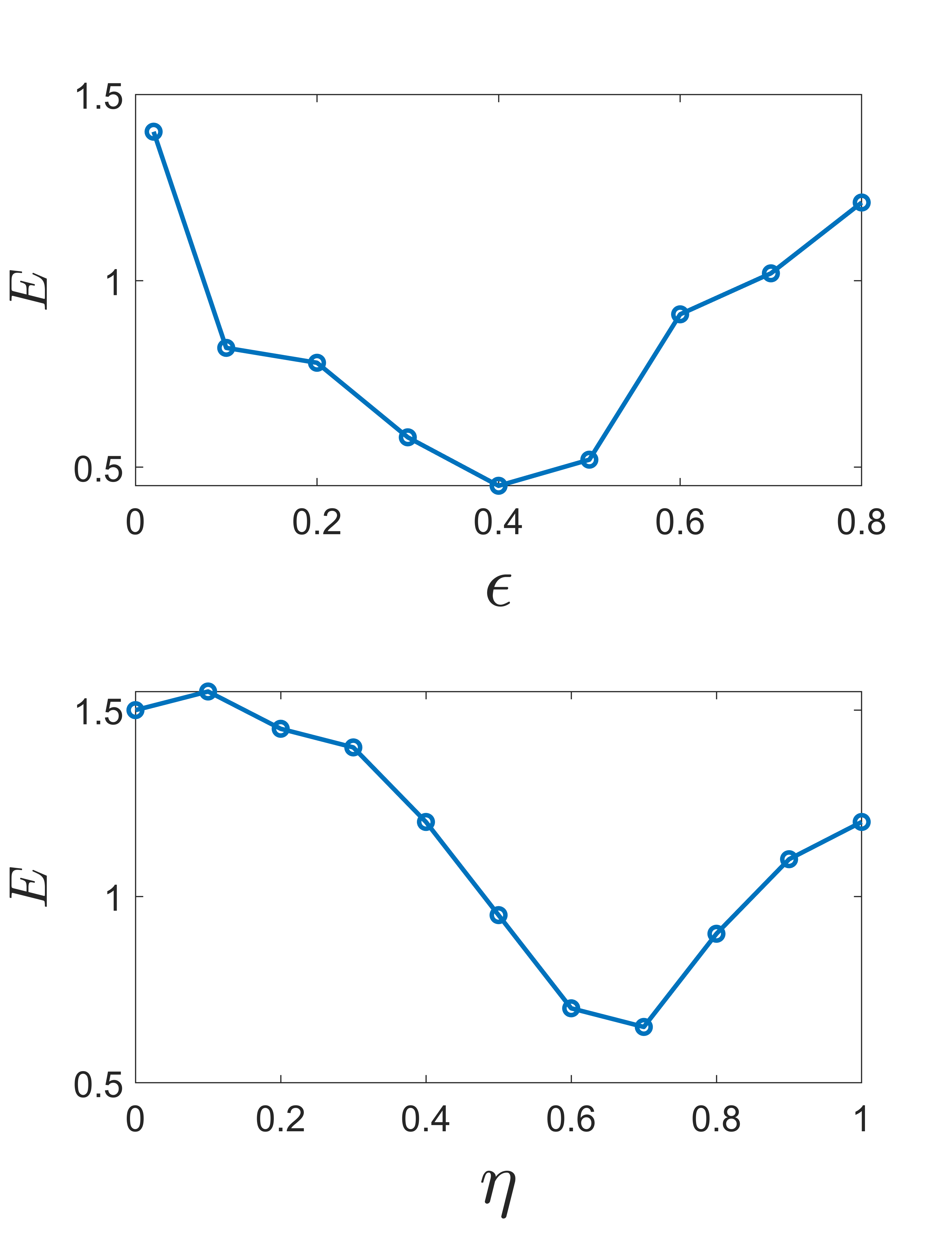}}
\end{subfigure}
	\caption{Parameter sensitivity analysis of $\epsilon$ and $\eta$ using the UAV model.}
	\label{fig:Exp2}
\end{figure*}

As mentioned in Sec.~\ref{sec:proposed}, the thresholds $\epsilon$ and $\eta$ are tuned based on the specific application and the desired performance level. This section examines how varying these thresholds affects the proposed approach's performance and computational efficiency. For different values of $\epsilon$ and $\eta$, we calculate the average cost ($E$), the average computational time, and the convergence rate. The results are illustrated in Fig.~\ref{fig:Exp2}. The top row shows the results for different values of $\epsilon$, and the bottom row presents those for $\eta$. In each row, we plot (a) the average computational time, (b) the convergence rate, and (c) the average cost $E$. This experiment is conducted using the UAV model.

As shown in the top row of Fig.~\ref{fig:Exp2} (a) and (b), increasing $\epsilon$ reduces the average computational time and, at the same time, increases the convergence rate. Increasing $\epsilon$ makes the genetic optimization treat any solution that costs less than $\epsilon$ as optimal. Thus, it terminates the optimization's iterations early and counts this termination as a converge to the optimal solution before the termination of the cycle's time. Therefore, the average computational time decreased, and the convergence rate increased, as shown in Fig.~\ref{fig:Exp2} (a) and (b).

For the average cost $E$, an increase in $\epsilon$ reduces $E$ initially, i.e., enhances the performance, as shown in Fig.~\ref{fig:Exp2} (c). However, there is a point at which further increasing $\epsilon$ causes diminishing returns, and the performance gets worse. This is because, in our proposed approach, we set the optimal search space to be around the previously found solution if it is optimal, i.e., its cost is less than $\epsilon$. Otherwise, the optimal search space is equal to the physical margin. Thus, increasing the value of $\epsilon$ enhances the probability of setting the margin to be around the previously found solution but leads to diminishing returns, as the proposed approach may set the margin around less precise optimal solutions, allowing errors in the control inputs to accumulate with the progression of control cycles. It is worth noting that, in some control scenarios, obtaining quick \textit{sub-optimal} control inputs is more critical than obtaining the best ones from the NMPC. These scenarios occur when running the NMPC while constrained with a short control cycle. In these cases, the threshold $\epsilon$ must be set to larger values to reduce computations, as adaptively setting $\BSM$  reduces the search space dramatically compared to the system's physical constraints.

On the other hand, increasing $\eta$ affects the average computational time differently than $\epsilon$, as shown in the bottom row of Fig.~\ref{fig:Exp2} (a). Specifically, increasing $\eta$ enhances the computational time initially because the margin used by the GA for sampling the solutions is set to the estimated BSM more often, in accordance with~\eqref {eq:final-bsm}. So, the number of required iterations is reduced as the GA searches within a smaller search space, decreasing the average computational time since the optimization process is usually terminated earlier. Nevertheless, increasing $\eta$ too much can eventually increase the computational cost. This happens because the GA may search for solutions within incorrect search space, becoming wrongly confident at the same time. As a result, computational time increases due to the computations spent searching in the wrong space, which is reflected in the convergence rate, as seen in Fig.~\ref{fig:Exp2} (b).

The average cost $E$ is also affected by the choice of $\eta$ as shown in Fig.~\ref{fig:Exp2} (c). At small values of $\eta$, the proposed approach behaves similarly to the OG due to the conservative estimation of the optimal BSM. However, as $\eta$ increases, the proposed approach becomes more confident in its estimations of the optimal BSM, leading to improved performance, and as a result, the average cost $E$ decreases. Nevertheless, increasing $\eta$ beyond a certain threshold causes diminishing returns and even negatively impacts performance, where the average cost $E$ increased. This occurs because a high $\eta$ value may cause the proposed approach to become overly reliant on predictions, even when wrong. In such cases, the proposed approach might trust wrong predictions, leading to incorrect decisions and a higher likelihood of accumulating errors over successive control cycles.

Thus, moderate values of $\eta$ balance trusting predictions and maintaining accuracy. Higher values reduce the robustness of the approach, as the proposed approach becomes too confident in predictions that may not be reliable. Therefore, $\eta$ must be chosen carefully to balance computational efficiency and prediction accuracy, avoiding overconfidence in less reliable predictions.
\section{Conclusion and future work}
\label{sec:conclusion}
This paper introduces a novel approach to accelerate genetic optimization in Nonlinear Model Predictive Control (NMPC) by dynamically learning the optimal search space size. The method trains a regression model to predict the most efficient search space in each control cycle, thereby increasing the likelihood of finding better optimal control inputs in the shortest computational time. We evaluated this approach on three nonlinear systems, comparing it with four other evolutionary optimization algorithms on an Nvidia\textsuperscript{\texttrademark} Jetson TX2 GPU in a processor-in-the-loop (PIL) setup. The results demonstrated that our approach outperformed the others, reducing computational time by [17-45\%]. It also significantly increased convergence to better control inputs within the cycle time by [35-47\%], yielding substantial performance improvements.

For future work, several promising directions emerge. One potential path is integrating the proposed method with the adaptive horizon estimation technique outlined in~\cite{FastRegMPC:Eslam:2022} to enhance NMPC's overall performance. Another avenue could be extending the proposed approach to other evolutionary optimization algorithms, allowing them to reduce their search spaces and achieve optimal solutions with faster computation times. Lastly, the proposed method could be applied to additional real-time industrial applications, where optimizing control inputs with both high accuracy and low computational cost is crucial.


\end{document}


\title{Supplementary Material for ``Accelerating genetic optimization of nonlinear model predictive control by learning optimal search space size''}
\author{Eslam~Mostafa,
        Hussein~A.~Aly,
        Ahmed~Elliethy
}

\maketitle
\IEEEpeerreviewmaketitle

\newcommand{\Figs}{../figs}
\renewcommand{\thefigure}{S.\arabic{figure}}
\renewcommand{\thetable}{S.\Roman{table}}
\renewcommand{\thesection}{S.\Roman{section}}
\newcommand{\grx}{\chi}
\newcommand{\gry}{\zeta}
\newcommand{\rgrx}{^r\grx}
\newcommand{\rgry}{^r\gry}
\newcommand{\betaB}{{\boldsymbol\beta}}
\newcommand{\omegaB}{{\boldsymbol\omega}}
\newcommand{\gammaB}{{\boldsymbol\gamma}}
\newcommand{\alphaB}{{\boldsymbol\alpha}}
\newcommand{\states}{\mathbf{x}}
\newcommand{\refx}{\mathbf{r}}
\newcommand{\refxi}{r}
\newcommand{\refu}{\mathbf{v}}
\newcommand{\refui}{v}
\newcommand{\optIns}{\mathbf{z}_c^*}
\newcommand{\ins}{\mathbf{z}_c}
\newcommand{\hzn}{h}
\newcommand{\bsm}{\delta_c}
\newcommand{\BSM}{\boldsymbol\Delta_{c}}
\newcommand{\BSMmax}{\Delta_{c}^{\max}}
\newcommand{\Cmax}{\mathcal{C}_{c}^{\max}}
\newcommand{\Emax}{\mathcal{E}_{c}^{\max}}
\newcommand{\Cc}{\mathbf{\mathcal{C}}_{c}}
\newcommand{\Ec}{\mathbf{\mathcal{E}}_{c}}
\newcommand{\Wc}{\mathbf{\mathcal{W}}_{c}}
\newcommand{\Rc}{\mathbf{\mathcal{R}}_{c}}
\newcommand{\kernelScale}{\gamma}
\newcommand{\SampDens}{\nu}

\newcommand{\ArxivFig}{./ArxivFig}

This supplementary document complements our main manuscript. The document provides additional details about the experimental setup and provides the mathematical formulations of the models used in the main manuscript. In Sec.~\ref{s.sec:UAV-mdl},Sec.~\ref{s.sec:VEC-mdl} and Sec.~\ref{s.sec:Rob-mdl}, we present the mathematical models of the Unmanned Aerial Vehicle (UAV), the ground vehicle, and the Single Link Flexible Joint Robot (SFJR), respectively, along with the values of their parameters. Finally, Section.~\ref{s.sec:sim-settings} outlines our proposed approach hyper-parameter settings used in all our experiments. 

\section{Unmanned aerial vehicle}
\label{s.sec:UAV-mdl}

We consider here an unmanned aerial vehicle (UAV) with four rotors~\cite{luukkonen2011UAVmdl}. Its motion is controlled by adjusting the angular velocities of these four rotors, i.e., the inputs are a vector of size $n=4$ that contains the angular velocities of the four rotors. The states of the UAV are determined by position, angular velocities, linear velocities, and orientations of the UAV, where each has three values corresponding to the 3D $X$, $Y$, and $Z$ dimensions. Thus, the state vector has size $m=12$ values~\cite{luukkonen2011UAVmdl}. The lose function of the UAV's model is constructed to penalize the error between the current UAV states with respect to their reference position states and penalize the velocities of the four rotors to save power. Additionally, the terminal cost $V$ is designed to be zero and we use using relatively long horizon approaches for stability insurance\cite{boccia2014longH_1,grimm2005longH_2,grune2010longH_3}. The loss function and the terminal cost at time instant $c$ are
\begin{equation}
\begin{aligned}
&\mathcal{L}(\bold{x}_c,\bold{z}_c)=\sum_{k=c}^{c+h-1}\left(\left[\refx_{k+1}-\mathbf{x}_{k+1}\right]^{T} \mathbf{Q}\left[\refx_{k+1}-\mathbf{x}_{k+1}\right]\right.\left.+\left[\refu_k-\mathbf{u }_{k}\right]^{T} \mathbf{R} \left[\refu_k-\mathbf{u }_{k}\right]\right), \\
&V = \mathbf{x}_{k+h}^{T} \Bar{\mathbf{Q}} \mathbf{x}_{k+h}
\end{aligned}
\label{eqn:MPCCost}
\end{equation}
where, $\refx_{k} = [\refxi^1_k, \dots, \refxi^{m}_{k}] ^{T}$ and $\refu_{k} = [\refui^1_k, \dots, \refui^{n}_{k}] ^{T}$ are the reference states and reference inputs at $k$, respectively. The matrix $\mathbf{Q}=\text{diag}(q^1,\dots,q^m)$ weighs the error between the reference and the system states according to the importance of each element of the state vector. The matrix $\mathbf{R}$ weighs the importance of the angular velocities of the four rotors. The matrix $\Bar{\mathbf{Q}}=\text{diag}(\Bar{q^1},\dots,\Bar{q^m})$ is the terminal cost weight.

\subsection{UAV mathematical model}
\begin{figure}[t!]
	\centerline{\includegraphics[scale=.35]{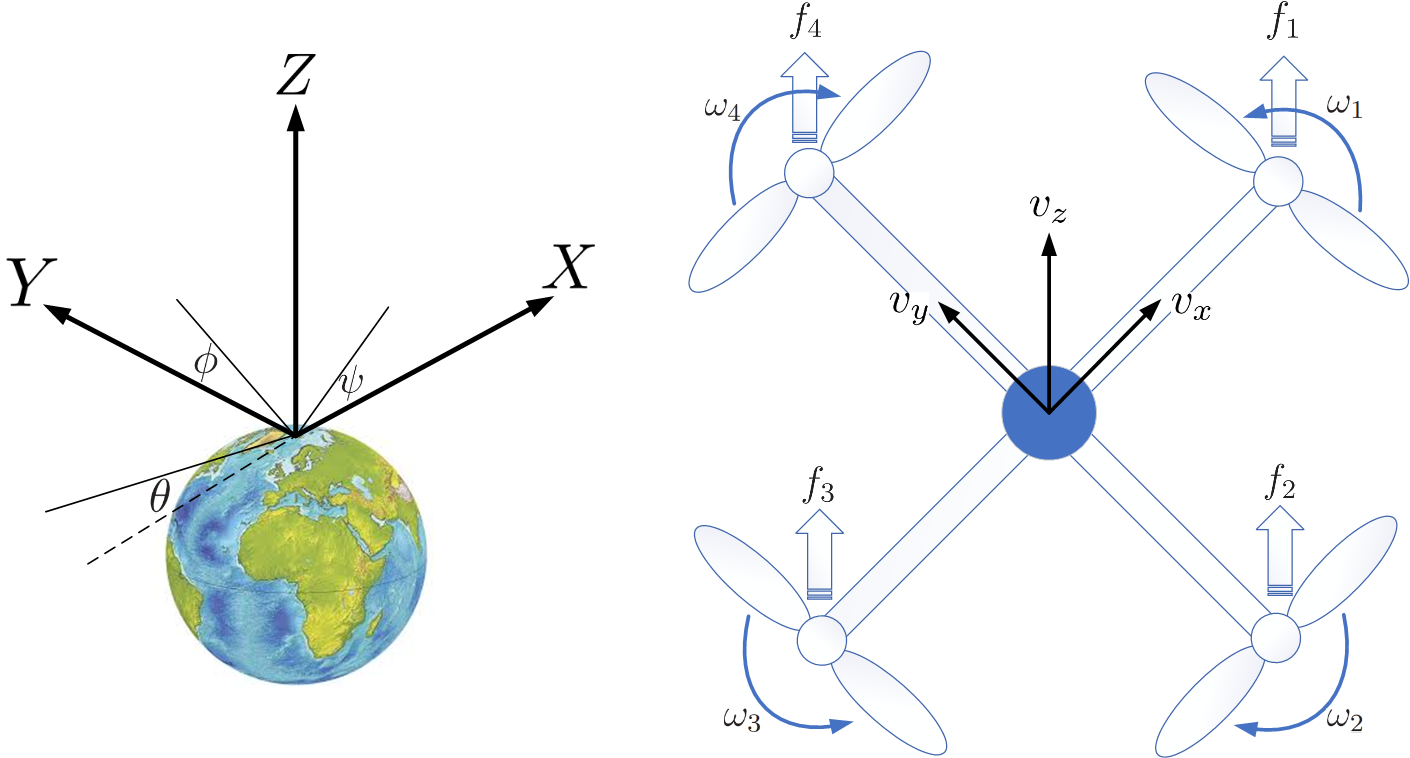}}
	\caption{Graphical representation of the reference frames, angular velocities, and forces generated by the four rotors in an UAV.}
	\label{fig:UAVmdl}
\end{figure}

Figure~\ref{fig:UAVmdl} describes a UAV in two reference frames: a fixed Body frame ($B$-frame) and a fixed earth reference frame ($E$-frame). The $B$-frame is attached to the body itself, and its origin is the center of mass of the UAV. The $E$-frame is a global reference frame with its origin at the center of the Earth and three orthogonal axes fixed to the Earth. UAV has six degrees of freedom, three describing its linear coordinates in the $X, Y,$ and $Z$-axes of the inertial frame by vector $\boldsymbol{\mu} = [X, Y, Z]^{T}$, and three describing its angular positions for inertial frame axes by vector $\boldsymbol{e} = [\phi, \vartheta, \psi]^{T}$. Roll angle $\phi$ determines the rotation around the $x$-axis. pitch angle $\vartheta$ around the $y$-axis and yaw angle $\psi$ around the $z$-axis. The linear velocities of the UAV are determined in the body frame by vector $\mathbf{v}_B = [v_x, v_y, v_z]^{T}$, and the angular velocities by vector $\mathbf{r}=[p, q, r]^{T}$.
The UAV dynamic model consists of twelve states $\mathbf{x}=[X, Y, Z, \dot{X}, \dot{Y}, \dot{Z}, \phi, \vartheta, \psi, \dot{\phi}, \dot{\vartheta}, \dot{\psi}]^T$, and four control inputs $\mathbf{u}=[\omega_1, \omega_2, \omega_3, \omega_4]^T$ which are the four angular velocity of its four rotors. The matrices $\mathbf{R}$ and $\mathbf{T}$ are the rotational matrix and translation matrix~\cite{metni2006rot-trans}, respectively, used to map the $B$-frame to the $E$-frame, where
\begin{equation*}
	\mathbf{R}=\left[\begin{array}{ccc}
		C_\psi C_\vartheta & C_\psi S_\vartheta S_\phi-S_\psi C_\phi & C_\psi S_\vartheta C_\phi+S_\psi S_\phi \\
		S_\psi C_\vartheta & S_\psi S_\vartheta S_\phi+C_\psi C_\phi & S_\psi S_\vartheta C_\phi-C_\psi S_\phi \\
		-S_\vartheta & C_\vartheta S_\phi & C_\vartheta C_\phi
	\end{array}\right],
\end{equation*}
\begin{equation*}
	\mathbf{T}=\left[\begin{array}{ccc}
		1 & 0 & -S_\vartheta \\
		0 & C_\phi & C_\vartheta S_\phi \\
		0 & -S_\phi & C_\vartheta C_\phi\\
	\end{array}\right],
\end{equation*}
$S_{a} = \sin(a)$, $C_{a} = \cos(a)$, and $a \in\{\phi, \psi, \vartheta\}$. The UAV is assumed to have symmetric structure with the four rotors aligned with the body $x$-axis and $y$-axis. Thus, the inertia matrix is a diagonal matrix
\begin{equation*}
	\mathbf{I} = \left[\begin{array}{ccc}
		I_{xx}  &0      &0  \\
		0       &I_{yy} &0 \\
		0       &0      &I_{zz}
	\end{array}\right],
\end{equation*}
in which $I_{xx} = I_{yy}$. The angular velocity and acceleration of the $i^{\text{th}}$ rotor create torque $\tau_{M_i}$ around the rotor axis that is defined as
\begin{equation*}
	\begin{aligned}
		\tau_{M_i}= b\omega_i^2 + I_M \dot{\omega},
	\end{aligned}
\end{equation*}
where $\omega_i$ and $\dot{\omega}_i$ is the angular velocity and acceleration, respectively of the $i^{\text{th}}$ rotor. $b$ is the drag constant and the $I_M$ is moment of inertia of all rotors. $\dot{\omega}_i$ can be neglected since it usually small value. Also, the angular velocity creates force $f_i = k\omega_i^2 $ in the direction of the rotor axis where $k$ is the lift constant. This forces create thrust $T_r$ in the direction of the body $z$-axis, which defined as
\begin{equation*}
	\begin{aligned}
		&T_r=&\sum_{i=1}^{4}f_i.
	\end{aligned}
\end{equation*}
Torque $\boldsymbol{\tau}$ consists of the torques $\tau_\phi, \tau_\vartheta$ and $\tau_\psi$ in the direction of the corresponding body frame angles, which defined as
\begin{equation*}
	\begin{aligned}
		&\boldsymbol{\tau} = &\left[ \begin{array}{c}
			\tau_\phi \\
			\tau_\vartheta \\
			\tau_\psi
		\end{array} \right]= & \left[\begin{array}{c}
			lk(-\omega_2^2 + \omega_4^2) \\
			lk(-\omega_1^2 + \omega_3^2) \\
			\tau_{M_i} 
		\end{array}\right],
	\end{aligned}
\end{equation*}
where $l$ is the distance between the rotor and the center of mass of the UAV.

To drive the nonlinear differential equations that describe the UAV model, we show the following relations. First, the relation between forces in the body frame is
\begin{equation*}
	m \dot{\mathbf{v}}_B+\mathbf{r} \times\left(m \mathbf{v}_B\right)=\mathbf{R}^{\mathrm{T}} \boldsymbol{G}+\mathbf{T}_B,
\end{equation*}
where $\mathbf{T}_B = [0, 0, T_r]^{T}$ is the total thrust of the rotors, $\boldsymbol{G} = [0,0,g]^T$, $g$ is the gravitational acceleration, and $m$ is the mass of the UAV. Also, the summation of the angular acceleration of the inertia $\mathbf{I} \dot{\mathbf{r}}$, the centripetal forces $\mathbf{r} \times(\mathbf{I} \mathbf{r})$, and the gyroscopic forces $\mathbf{\Gamma}$ are equal to the external torque $\boldsymbol{\tau}$, i.e.,
\begin{equation}
	\mathbf{I} \dot{\mathbf{r}}+\mathbf{r} \times(\mathbf{I} \mathbf{r})+\mathbf{\Gamma}=\boldsymbol{\tau}.
\end{equation}
Therefore,
\begin{equation}
	\begin{gathered} \\
		\dot{\mathbf{r}}=\mathbf{I}^{-1}\left(-\left[\begin{array}{c}
			p \\
			q \\
			r
		\end{array}\right] \times\left[\begin{array}{c}
			I_{x x} p \\
			I_{y y} q \\
			I_{z z} r
		\end{array}\right]-I_r\left[\begin{array}{c}
			p \\
			q \\
			r
		\end{array}\right] \times\left[\begin{array}{l}
			0 \\
			0 \\
			1
		\end{array}\right] \omega_{\Gamma}+\boldsymbol{\tau}\right), \\
		{\left[\begin{array}{c}
				\dot{p} \\
				\dot{q} \\
				\dot{r}
			\end{array}\right]=\left[\begin{array}{l}
				\left(I_{y y}-I_{z z}\right) q r / I_{x x} \\
				\left(I_{z z}-I_{x x}\right) p r / I_{y y} \\
				\left(I_{x x}-I_{y y}\right) p q / I_{z z}
			\end{array}\right]-I_r\left[\begin{array}{c}
				q / I_{x x} \\
				-p / I_{y y} \\
				0
			\end{array}\right] \omega_{\Gamma}+\left[\begin{array}{c}
				\tau_\phi / I_{x x} \\
				\tau_\vartheta / I_{y y} \\
				\tau_\psi / I_{z z}
			\end{array}\right]},
	\end{gathered}
	\label{eq:bodyEqn1}
\end{equation}
where $\omega_{\Gamma}=\omega_1-\omega_2+\omega_3-\omega_4$. The angular accelerations in the inertial frame are then attracted from the body frame accelerations with the transformation matrix $\mathbf{T}^{-1}$ and its time derivative as
\begin{equation}
	\begin{aligned}
		\ddot{\boldsymbol{e}} & =\frac{\mathrm{d}}{\mathrm{d} t}\left(\mathbf{T}^{-1} \mathbf{r}\right)=\frac{\mathrm{d}}{\mathrm{d} t}\left(\mathbf{T}^{-1}\right) \mathbf{r}+\mathbf{T}^{-1} \dot{\mathbf{r}}, \\
		= & {\left[\begin{array}{ccc}
				0 & \dot{\phi} C_\phi T_\vartheta+\dot{\vartheta} S_\phi / C_\vartheta^2 & -\dot{\phi} S_\phi C_\vartheta+\dot{\vartheta} C_\phi / C_\vartheta^2 \\
				0 & -\dot{\phi} S_\phi & -\dot{\phi} C_\phi \\
				0 & \dot{\phi} C_\phi / C_\vartheta+\dot{\phi} S_\phi T_\vartheta / C_\vartheta & -\dot{\phi} S_\phi / C_\vartheta+\dot{\vartheta} C_\phi T_\vartheta / C_\vartheta
			\end{array}\right] \mathbf{r}+\mathbf{T}^{-1} \dot{\mathbf{r}} },
	\end{aligned}
	\label{eq:bodyEqn}
\end{equation}
where the angular velocities are transformed from $B$-frame to the $E$-frame using the translation matrix $\mathbf{T}$.

In the inertial frame, only the gravitational force and the magnitude and direction of the thrust affects the acceleration of the UAV because the centrifugal force is repealed, i.e., 
\begin{equation*}
	m \ddot{\boldsymbol{p}}=\boldsymbol{G}+\mathbf{r} \boldsymbol{T}_{B}.
\end{equation*}
Thus,
\begin{equation}
	\begin{aligned}
		{\left[\begin{array}{c}
				\ddot{X} \\
				\ddot{Y} \\
				\ddot{Z}
			\end{array}\right]=
			-g\left[\begin{array}{l}
				0 \\
				0 \\
				1
			\end{array}\right]+
			\frac{T_r}{m}\left[\begin{array}{c}
				C_\psi S_\vartheta C_\phi+S_\psi S_\phi \\
				S_\psi S_\vartheta C_\phi-C_\psi S_\phi \\
				C_\vartheta C_\phi
			\end{array}\right] .}
	\end{aligned}
\end{equation}
As in reality, drag force generated by the air resistance should be included by including the diagonal coefficient matrix associating the linear velocities to the force which slowing the movement, i.e.,
\begin{equation}
	\begin{aligned}
		\left[\begin{array}{c}
			\ddot{X} \\
			\ddot{Y} \\
			\ddot{Z}
		\end{array}\right]=-g\left[\begin{array}{l}
			0 \\
			0 \\
			1
		\end{array}\right]+\frac{T_r}{m}\left[\begin{array}{c}
			C_\psi S_\vartheta C_\phi+S_\psi S_\phi \\
			S_\psi S_\vartheta C_\phi-C_\psi S_\phi \\
			C_\vartheta C_\phi
		\end{array}\right]-\frac{1}{m}\left[\begin{array}{ccc}
			A_x & 0 & 0 \\
			0 & A_y & 0 \\
			0 & 0 & A_z
		\end{array}\right]\left[\begin{array}{l}
			\dot{X} \\
			\dot{Y} \\
			\dot{Z}
		\end{array}\right] ,
	\end{aligned}
	\label{eq:InertialEqn}
\end{equation}
where $A_x$, $A_y$ and $A_z$ are the drag force coefficients for velocities in the corresponding directions of the inertial frame. From the above discussion, the nonlinear differential equations that describe the UAV model are the equations~\ref{eq:bodyEqn} and ~\ref{eq:InertialEqn}. The parameter values of the UAV model are listed in Table~\ref{tab:UAVparameters}. Finally, the UAV hyper-parameters used in our NMPC optimization are presented in Table~\ref{tab:UAV-parms}.
\begin{table}[]
	\caption{Parameter values of the UAV model}
	\centering
	\setlength{\tabcolsep}{2mm}
	\renewcommand{\arraystretch}{1.5}
	\setlength{\arrayrulewidth}{0.2mm}
	\begin{tabular}{|c|c|c|}
		\hline Parameter & Value & Unit \\
		\hline$g$ & $9.81$ & $\mathrm{~m} / \mathrm{s}^2$ \\\hline
		$m$ & $0.468$ & $\mathrm{~kg}$ \\\hline
		$l$ & $0.225$ & $\mathrm{~m}$ \\\hline
		$k$ & $2.980 \cdot 10^{-6}$ & -\\\hline
		$b$ & $1.140 \cdot 10^{-7}$ & -\\\hline
		$I_M$ & $3.357 \cdot 10^{-5}$ & $\mathrm{~kg} \mathrm{~m}^2$ \\\hline
		
		$I_{x x}$ & $4.856 \cdot 10^{-3}$ & $\mathrm{~kg} \mathrm{~m}^2$ \\\hline
		$I_{y y}$ & $4.856 \cdot 10^{-3}$ & $\mathrm{~kg} \mathrm{~m}^2$ \\\hline
		$I_{z z}$ & $8.801 \cdot 10^{-3}$ & $\mathrm{~kg} \mathrm{~m}{ }^2$ \\\hline
		$A_x$ & $0.25$ & $\mathrm{~kg} / \mathrm{s}$ \\\hline
		$A_y$ & $0.25$ & $\mathrm{~kg} / \mathrm{s}$ \\\hline
		$A_z$ & $0.25$ & $\mathrm{~kg} / \mathrm{s}$ \\\hline
	\end{tabular}
	\label{tab:UAVparameters}
\end{table}

\begin{table}
	\centering
	\caption{UAV hyper-parameters used in the NMPC optimization}
	\setlength{\tabcolsep}{2mm}
	\renewcommand{\arraystretch}{1.5}
	\setlength{\arrayrulewidth}{0.2mm}
	\centering
	\begin{tabular}{|c|c|c|}
		\hline Parameter & value     \\ \hline
		$t_s$ & $.02s$ \\\hline 
		$h$ & $10$   \\\hline
		$\mathbf{Q}$ &  diag(1, 1, 1, 0, 0, 0, 1, 1, 1, 0, 0, 0)  \\\hline
		$\mathbf{R}$ & diag(.1, .1, .1, .1) \\\hline
		
		$\mathbf{u_{\text{min} } }$ & [0,0,0,0] \\ \hline
		$\mathbf{u_{\text{max} } }$ & [12,12,12,12]\\ \hline
		$\mathbf{\Delta u_{\text{min}} }$ & [-.2,-.2] \\ \hline
		$\mathbf{ \Delta u_{\text{max}}} $ & [.2,.2]  \\\hline
		
		\multirow{2}{*}{$\mathbf{x_{\text{min} } }$}   
		& $[-\inf, -\inf, -\inf, -\inf, -\inf, -\inf,$ \\
		& $-\frac{\pi}{3}, -\frac{\pi}{3}, -\frac{\pi}{3},  -\frac{\pi}{24}, -\frac{\pi}{24}, -\frac{\pi}{24}]  $ \\ \hline
		$\mathbf{x_{\text{max} } } $ 
		& $[\inf, \inf, \inf, \inf, \inf, \inf,  \frac{\pi}{3}, \frac{\pi}{3}, \frac{\pi}{3}, \frac{\pi}{24}, \frac{\pi}{24}, \frac{\pi}{24}]  $ \\ \hline		
	\end{tabular}
	\label{tab:UAV-parms}
\end{table}

\section{Ground vehicle}
\label{s.sec:VEC-mdl}
Here, we consider an autonomous vehicle, where its motion is controlled by adjusting its steering angle and acceleration, i.e., the inputs are a vector of size $n=2$. Longitudinal and lateral positions, longitudinal and lateral linear velocities, yaw angle, and yaw rate determine the states of the vehicle. Thus, the state vector has size $m=6$ values~\cite{kong2015VECmdl}. The lose function of the NMPC for this model is designed to penalize any discrepancies between the vehicle's longitudinal and lateral positions and orientation and their corresponding references. Moreover, it ensures the smoothness of control input changes by adhering to imposed constraints on both states and inputs~\cite{du2016develGAsaftyEqn}. The terminal cost is formulated as predefined in~\eqref{eqn:MPCCost}. However, the lose function is formulated as predefined in~\eqref{eqn:MPCCost} with a minor modification to account for the safety and comfort of human passengers. Specifically, when the NMPC algorithm encounters certain conditions, such as emergency maneuvers or sudden changes in driving conditions, it may be necessary to adjust the weight matrices $\mathbf{Q}$ and $\mathbf{R}$ in the lose function~\eqref{eqn:MPCCost} to achieve the desired performance. Let us define $cond1$ to be true whenever any of these conditions occur, then the $\mathbf{Q}$ and $\mathbf{R}$ matrices are adjusted as
\begin{equation}
\begin{aligned}
&\mathbf{Q}= \begin{cases}
             \mathbf{Q}_1              & \text{if } cond1,\\
             \mathbf{Q}_0             & \text{otherwise},
        \end{cases} \\
&\mathbf{R}= \begin{cases}
             \mathbf{R}_1              & \text{if } cond1,\\
             \mathbf{R}_0              & \text{otherwise}.
        \end{cases} \\
\end{aligned}
\label{eqn:MPCCost_update}
\end{equation}
Thus, based on $cond1$, the NMPC adjusts $\mathbf{Q}$ by selecting between  $\mathbf{Q}_0$ and $\mathbf{Q}_1$, where each matrix differently weighs the error between the reference and the system states according to the importance of each element of the state vector. Similarly, $\mathbf{R}$ is adjusted. 

\begin{figure}[]
	\centerline{\includegraphics[scale=.3]{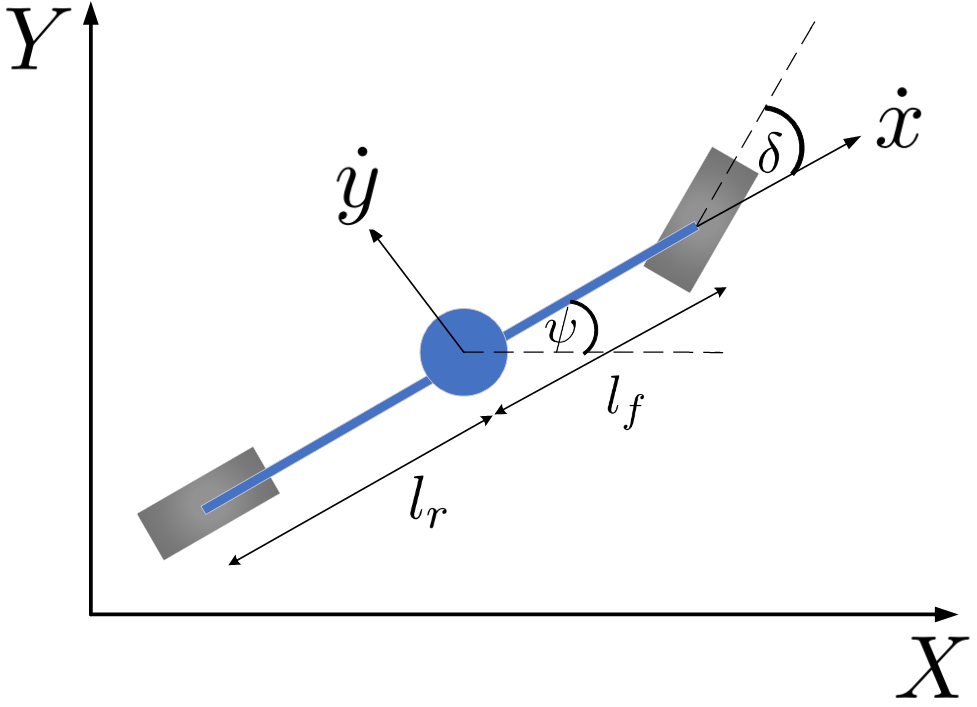}}
	\caption{A graphical representation of a dynamic bicycle model}
	\label{fig:VECmdl}
\end{figure}

Vehicle system control involves the critical task of regulating and maintaining the longitudinal position, lateral position, and orientation of the vehicle. These three factors are crucial for ensuring the safe and efficient operation of the vehicle. We use a dynamic bicycle model to represent the vehicle dynamics~\cite{kong2015VECmdl}. The dynamic bicycle model is a simplified vehicle model commonly used to analyze and control vehicle dynamics. It provides a good balance between accuracy and simplicity. The nonlinear differential equations described by the dynamic bicycle model are
\begin{equation}
	\begin{aligned}
		& \dot{X}=\dot{x} \cos \psi-\dot{y} \sin \psi \\
		& \dot{Y}=\dot{x} \sin \psi+\dot{y} \cos \psi \\
		& \ddot{x}=\dot{y}+a \\
		& \ddot{y}=-\dot{\psi} \dot{x}+\frac{2}{m}\left(F_{c, f} \cos \delta+F_{c, r}\right) \\
		& \ddot{\psi}=\frac{2}{I_z}\left(l_f F_{c, f}-l_r F_{c, r}\right) , \\
	\end{aligned}
\end{equation}
where $\dot{x}$ and $\dot{y}$ denote the longitudinal and lateral speeds in the body frame, respectively, and $\dot{\psi}$ denotes the yaw rate. $m$ and $I_z$ denote the vehicle's mass and yaw inertia, respectively. $F_{c, f}$ and $F_{c, r}$ denote the lateral tire forces at the front and rear wheels, respectively. For the linear tire model, $F_{c, i}$ is defined as $$F_{c, i}=-C_{\alpha_i} \alpha_{i} ,$$ where $i \in\{f, r\}$, $\alpha_i$ is the tire slip angle, and $C_{\alpha_i}$ is the tire cornering stiffness. A graphical representation of a dynamic bicycle model is shown in Figure~\ref{fig:VECmdl}.

This model consists of six states $\mathbf{x}=[X, Y, \dot{x}, \dot{y}, \dot{\psi}, \psi]^T$, which describe the vehicle's longitudinal and lateral positions with respect to the inertial frame, longitudinal and lateral velocities, yaw angle and yaw rate. The system employs two control inputs: steering angle and acceleration $\mathbf{u}=[a, \delta]^T$. Table~\ref{tab:VECparameters} lists the vehicle's model parameter values. Finally, the Vehicle hyper-parameters used in our NMPC optimization are presented in Table~\ref{tab:VEC-parms}.
\begin{table}[]
	\centering
	\setlength{\tabcolsep}{2mm}
	\renewcommand{\arraystretch}{1.5}
	\setlength{\arrayrulewidth}{0.2mm}
	\caption{Vehicle model Parameters}
	\label{tab:VECparameters}
	\begin{tabular}{|c|c|c|c|}
		\hline Parameter            & Symbol         & Value             & Units \\
		\hline
		Vehicle mass  & ${m}$          & 1,650             & kg \\
		\hline Yaw inertia          & ${I}_{z}$      & 2,650     & kg$\cdot$ $\textnormal{m}^{2}$ \\
		\hline Front axle to CG     & ${l}_{f}$      & 1.1               & m \\
		\hline Rear axle to CG      & ${l}_{r}$      & 1.7               & m \\
		\hline Cornering stiffness of front-axle     & $C_f$ & 55,494    & N/rad \\
		\hline Cornering stiffness of rear-axle      & $C_{r}$ & 55,494   & N/rad \\
		\hline
	\end{tabular}
\end{table}

\begin{table}
	\centering
	\caption{Vehicle hyper-parameters used in NMPC optimization}
	\setlength{\tabcolsep}{2mm}
	\renewcommand{\arraystretch}{1.5}
	\setlength{\arrayrulewidth}{0.2mm}
	\centering
	\begin{tabular}{|c|c|c|}
		\hline Parameter & value     \\ \hline
		$t_s$ & $.02s$ \\\hline 
		$h$ & $10$   \\\hline
		$\mathbf{Q}$ &  diag(1, 1, 1, 0, 0, 0, 1, 1, 1, 0, 0, 0)  \\\hline
		$\mathbf{R}$ & diag(.1, .1, .1, .1) \\\hline
		
		$\mathbf{u_{\text{min} } }$ & [0,0,0,0] \\ \hline
		$\mathbf{u_{\text{max} } }$ & [12,12,12,12]\\ \hline
		$\mathbf{\Delta u_{\text{min}} }$ & [-.2,-.2] \\ \hline
		$\mathbf{ \Delta u_{\text{max}}} $ & [.2,.2]  \\\hline
		
		\multirow{2}{*}{$\mathbf{x_{\text{min} } }$}   
		& $[-\inf, -\inf, -\inf, -\inf, -\inf, -\inf,$ \\
		& $-\frac{\pi}{3}, -\frac{\pi}{3}, -\frac{\pi}{3},  -\frac{\pi}{24}, -\frac{\pi}{24}, -\frac{\pi}{24}]  $ \\ \hline
		$\mathbf{x_{\text{max} } } $ 
		& $[\inf, \inf, \inf, \inf, \inf, \inf,  \frac{\pi}{3}, \frac{\pi}{3}, \frac{\pi}{3}, \frac{\pi}{24}, \frac{\pi}{24}, \frac{\pi}{24}]  $ \\ \hline
		
  \end{tabular}
	\label{tab:VEC-parms}
\end{table}

\section{Single-link flexible-joint robot system}
\label{s.sec:Rob-mdl}

In this section, we describe the dynamics and control of the single-link flexible-joint robot (SFJR) system. This robot, operating in a vertical plane, is actuated by a DC motor with a gear reduction system. The system consists of two parts: the motor side and the link side, which are connected through an elastic element modeled as a linear spring. The motor side includes a DC motor and its drive, while the link side comprises a massless link and a load. The elastic coupling between the two sides introduces compliance in the system, modeled as a spring with stiffness $K$, which allows for joint deformation under torque~\cite{zhang2021DEcompare}. The angular position of the link side is controlled by adjusting the DC motor voltage, resulting in a single control input of size $n=1$. The state of the system includes the angular positions and angular velocities of both the motor and the link sides, along with the motor current, forming a state vector of size $m=5$. 
The lose function of the NMPC for this model is designed to penalize any discrepancies between the vehicle's longitudinal and lateral positions and orientation and their corresponding references. The terminal cost is formulated as predefined in~\eqref{eqn:MPCCost}.

The dynamics of the SFJR can be expressed using the Euler-Lagrange formulation. The angular positions of the motor and the link are denoted as $\varepsilon_2$ and $\varepsilon_1$, respectively, while the deformation of the spring is given by $\varphi = \varepsilon_2 - \varepsilon_1$. The governing equations of motion for the system are described as:
\begin{equation}
\begin{aligned}
    &J_1 \ddot{\varepsilon}_1 + mgl \sin(\varepsilon_1) + K_f1 \dot{\varepsilon}_1 = K (\varepsilon_2 - \varepsilon_1), \\
    &J_2 \ddot{\varepsilon}_2 + K_f2 \dot{\varepsilon}_2 + K (\varepsilon_2 - \varepsilon_1) = N K_\tau i, \\
    &R_m i + L \dot{i} + N K_e \dot{\varepsilon}_2 = U_v,
\end{aligned}
\end{equation}
where $J_1$ and $J_2$ are the moments of inertia of the link and motor side, $K_f1$ and $K_f2$ represent the viscous damping coefficients, and $N$ is the gear reduction ratio. The terms $K_\tau$, $K_e$, and $R_m$ refer to the motor torque constant, back electromotive force coefficient, and resistance of the motor's armature, respectively. The motor voltage is represented by $U_v$, while the current is denoted by $i$.

We define the state vector $\mathbf{x} = [\varepsilon_1, \dot{\varepsilon}_1, \varepsilon_2, \dot{\vartheta}_2, i]^T$ and control input $u = U_v$. The dynamic behavior of the robot is primarily influenced by the elasticity of the joint, which introduces challenges such as overshoots and residual vibrations during control~\cite{zhang2021DEcompare}.

The parameters used for the simulation of the SFJR system are listed in Table~\ref{tab:FJparameters}. Finally, the SFJR hyper-parameters used in our NMPC optimization are presented in Table~\ref{tab:SFJR-parms}.

\begin{table}[h]
    \centering
    \caption{Simulation model parameters of the SFJR.}
    	\setlength{\tabcolsep}{2mm}
	\renewcommand{\arraystretch}{1.5}
	\setlength{\arrayrulewidth}{0.2mm}
	\centering
    \label{tab:FJparameters}
    \begin{tabular}{|c|c|c|}
        \hline
        \textbf{Parameter} & \textbf{Symbol} & \textbf{Value} \\
        \hline
        Moment of inertia of the link side & $J_1$ & 0.8 kg$\cdot$m$^2$ \\        \hline

        Moment of inertia of the motor side & $J_2$ & 0.1 kg$\cdot$m$^2$ \\        \hline

        Viscous damping on link side & $K_{f1}$ & 2.0 \\        \hline

        Viscous damping on motor side & $K_{f2}$ & 2.0 \\        \hline

        Spring stiffness & $K$ & 70 Nm/rad \\        \hline

        Motor torque constant & $K_\tau$ & 9.3$\times 10^{-3}$ Nm/A \\        \hline

        Armature resistance & $R_m$ & 5.3 $\Omega$ \\        \hline

        Armature inductance & $L$ & 1.4$\times 10^{-5}$ H \\        \hline

        Back electromotive force constant & $K_e$ & 0.1 V/rad/s \\        \hline

        Gear ratio & $N$ & 200 \\        \hline

        Load mass & $m$ & 0.3 kg \\        \hline

        Link length & $l$ & 0.5 m \\        \hline

        Gravitational acceleration & $g$ & 9.8 m/s$^2$ \\
        \hline
    \end{tabular}
\end{table}

\begin{table}
	\centering
	\caption{SFJR hyper-parameters used in NMPC optimization}
	\setlength{\tabcolsep}{2mm}
	\renewcommand{\arraystretch}{1.5}
	\setlength{\arrayrulewidth}{0.2mm}
	\centering
	\begin{tabular}{|c|c|c|}
		\hline Parameter & value     \\ \hline
		$t_s$ & $.04s$  \\\hline
		$h$ & $10$  \\\hline
		$\mathbf{Q}$ &  diag(1, 0, 0, 0, 0)  \\\hline
		$\mathbf{R}$ & .5 \\\hline
		$\mathbf{u_{\text{min} } }$  &$0$\\ \hline
		$\mathbf{u_{\text{max} } }$  &$24$\\ \hline
		$\mathbf{\Delta u_{\text{min}} }$  & $ -.1$\\\hline
		$\mathbf{ \Delta u_{\text{max}}} $ & $ .1$\\\hline
		
		$\mathbf{x_{\text{min} } }$ & $[-\pi, \frac{-\pi}{18}, -\pi, \frac{-\pi}{18}, 0] $ \\\hline
		$\mathbf{x_{\text{max} } } $  & $[\pi, \frac{\pi}{18}, \pi, \frac{\pi}{18}, 5] $\\\hline
		
	\end{tabular}
	\label{tab:SFJR-parms}
\end{table}

\section{Experimental settings}
\label{s.sec:sim-settings}
\begin{figure}[]
	\centering
	\begin{subfigure}[]{.8\textwidth}
		\centerline{\includegraphics[width=\textwidth]{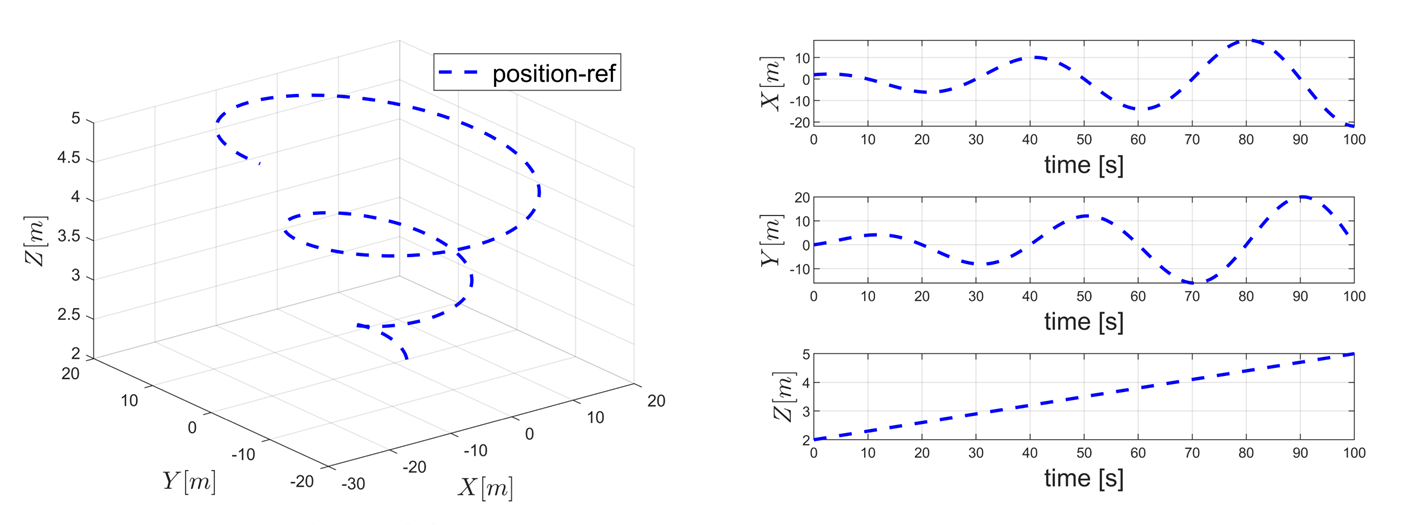}}
		\caption{}
	\end{subfigure}
	\begin{subfigure}[]{.8\textwidth}
	\centerline{\includegraphics[width=\textwidth]{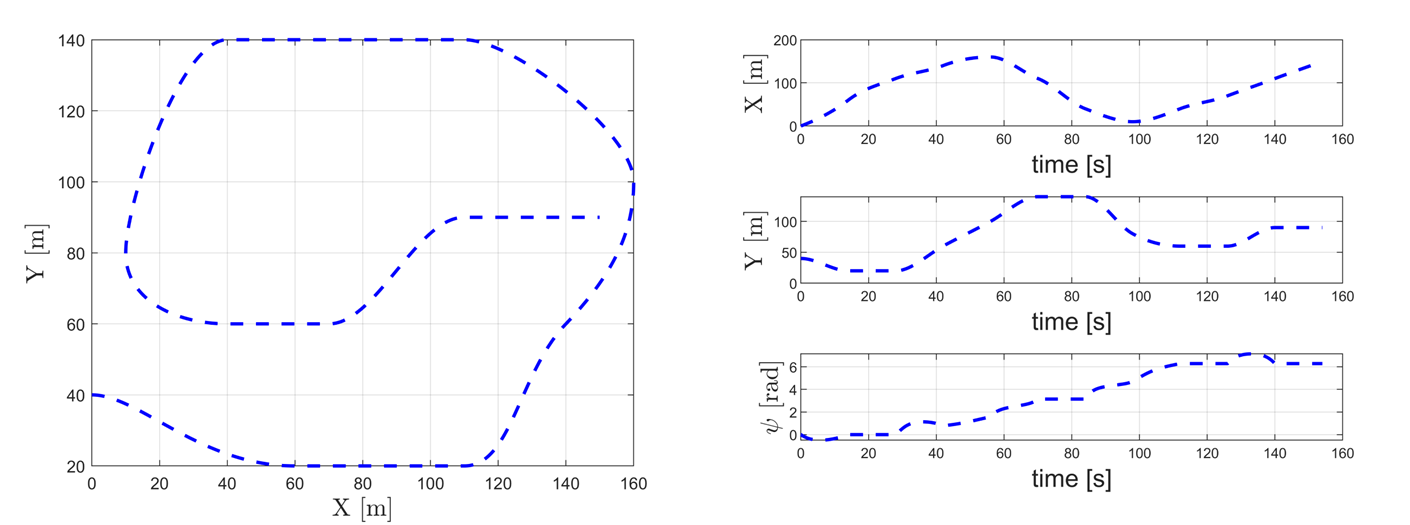}}
		\caption{}
	\end{subfigure}
 \begin{subfigure}[]{.4\textwidth}
	\centerline{\includegraphics[width=\textwidth]{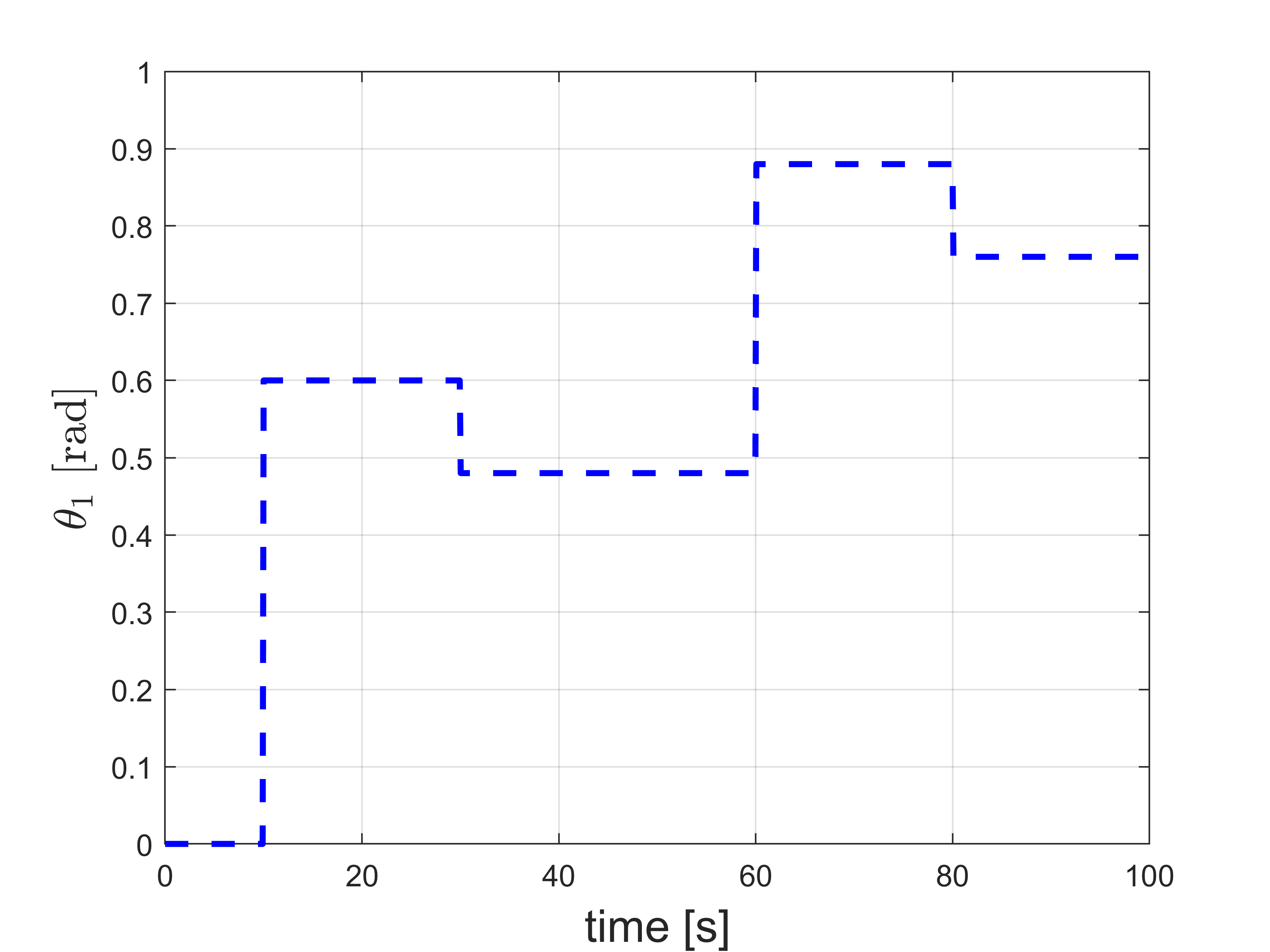}}
		\caption{}
	\end{subfigure}
	\caption{Reference states for the three test cases. UAV's reference states are shown in (a). Vehicle's reference states are shown in (b). Reference states for the SFJR are shown in (c).}
	\label{fig:s.refStates}
\end{figure}
This section presents the experimental settings employed in all experiments in the main menu script. The hyper-parameters of our regression model are estimated using 5-fold cross-validation~\cite{chorowski2014crossvalidation}. The training uses the dual problem form~\cite{vapnik1999dualformula}, and a Gaussian kernel is applied. 
The proposed approach hyper-parameter settings for the UAV, the ground vehicle and the SFJR are presented in Table~\ref{tab:parms}. The reference states for both the UAV, the vehicle and SFJR models are shown in Fig.~\ref{fig:s.refStates}. The simulation lengths for the UAV, vehicle models and the SFJR are set to 100, 100 and 160 seconds, respectively.
\begin{table}
	\centering
	\caption{Proposed approach hyper-parameter settings}
	\setlength{\tabcolsep}{2mm}
	\renewcommand{\arraystretch}{1.5}
	\setlength{\arrayrulewidth}{0.2mm}
	\centering
	\begin{tabular}{|c|c|c|c|}
		\hline Parameter & UAV& Vehicle & SFJR     \\ \hline
    	$\upsilon$ & $.95$ & .95 & .95 \\\hline
		$\epsilon$  & .4 & 2 & .5\\\hline
            $\eta$  & .7 & .65 & .8\\\hline
            $\rho$  & .1 & .1 & .2\\\hline
		$\theta$  & .05 & .15 & .15\\\hline
		$\SampDens$ & 100 & 200 & 200\\\hline 
		$\xi$ & 10 &  20 & 20\\\hline 
		$G$ & 3 & 5 & 4\\\hline
		Selection method &  Tournament selection & Tournament selection 
            & Tournament selection\\\hline
		Crossover rate & .4 & .5 & .4\\\hline
		Mutation rate & .05 & .1 & .1\\\hline
		$\lambda$ & [0.2,0.2,0.2,0.2]    & [0.1,0.35] & 0.5\\\hline
            $\kernelScale$ & [0.1,0.2,0.2,0.05]   & [0.05,.01] & 0.2\\\hline
            $C$ & [10,5,5,10]  & [1,2]  & 5         \\\hline
	\end{tabular}
	\label{tab:parms}
\end{table}